%On the integral cohomology of wreath products, by Ian J. Leary.  
%THIS IS A PLAIN TEX FILE.  
\magnification 1200
\input amssym.def
\input amssym.tex
\input xy 
\xyoption{v2} 

\def\wreath{\wr}

\def\proof{\medbreak \noindent{\bf Proof. }}
\def\remarks{\medbreak \noindent{\bf Remarks. }}
\def\remark{\medbreak \noindent{\bf Remark. }}

\def\mapright#1{\smash{\mathop{\longrightarrow}\limits^{#1}}}
\def\mapleft#1{\smash{\mathop{\longleftarrow}\limits^{#1}}}
\def\mapdown#1{\Big\downarrow
    \rlap{$\vcenter{\hbox{$\scriptstyle#1$}}$}}
\def\mapup#1{\Big\uparrow
    \rlap{$\vcenter{\hbox{$\scriptstyle#1$}}$}}
\def\mkpt#1\par{\item{$\bullet$} #1}

\def\beginsection#1\par{\vskip0pt plus.1\vsize\penalty-250
\vskip0pt plus-.1\vsize\bigskip\vskip\parskip
\message{#1}
\centerline{#1}\nobreak\smallskip\noindent}

\def\wre#1{{#1\wreath 1}}

\def\res {{\rm Res}}

\def\ind {{\rm Ind}}
\def\norm {{\cal N}}
\def\spot {\bullet}
\def\im {{\rm Im}}
\def\rk {{\rm Rk}}
\def\Hom {{\rm Hom}}
\def\Ext {{\rm Ext}}
\def\Tor {{\rm Tor}}
\def\tot {{\rm Tot}}
\def\zz {{\Bbb Z}}
\def\fp {{\Bbb F}_p}
\def\ftwo{{\Bbb F}_2}
\def\plocal {{\Bbb Z}_{(p)}}
\def\ece{{\rm ee}}
\def\ce{{\rm e}}
\def\eve{{\rm e}}
\def\odd{{\rm o}}
\def\pra{\par}
\def\qed{\hfill${\scriptstyle\blacksquare}$}
\def\inj{\rightarrowtail}
\def\surj{\twoheadrightarrow}

\def\Ann{{\rm Ann}}
\def\bu{\bullet}
\def\ci{\circ}
\def\mi{\cdot}
\def\blocksquare{\vrule height .9ex width .8ex depth -.1ex}
\def\sq{\blocksquare}
\newbox\bigstrutbox
\setbox\bigstrutbox=\hbox{\vrule height9.5pt depth3.5pt width0pt}
\def\bigstrut{\relax\ifmmode\copy\bigstrutbox\else\unhcopy\bigstrutbox\fi}
\def\ze{{\Bbb Z}}
\def\mj{\hbox{-}}

\def\book#1/#2/#3/#4/#5/#6/ {\item{[#1]} #2{\sevenrm #3}. {\it #4}.\
#5, (#6).\par\smallskip}
\def\paper#1/#2/#3/#4/#5/#6/#7/#8/ {\item{[#1]} #2{\sevenrm #3}. #4.
{\it #5} {\bf #6} (#7), #8.\par\smallskip}
\def\prepaper#1/#2/#3/#4/#5/ {\item{[#1]} #2{\sevenrm #3}. #4.
#5\par\smallskip}

\def\ade{1}
\def\adca{2}
\def\ademil{3}
\def\atimac{4}
\def\ben{5}
\def\beninv{6}
\def\bencarl{7}
\def\bro{8}
\def\brobui{9}
\def\lmm{10}
\def\carl{11}
\def\caei{12}
\def\eva{13}
\def\evb{14}
\def\evk{15}
\def\hako{16}
\def\hwh{17}
\def\huea{18}
\def\jrh{19}
\def\kah{20}
\def\ian{21}
\def\oprob{22}
\def\nak{23}
\def\pal{24}
\def\qui{25}
\def\quib{26}
\def\quies{27}
\def\spa{28}
\def\ste{29}

\font \smallfont=cmr8 at 8pt
\font \smallbold=cmbx8 at 8pt

\phantom{}
\vskip 1 truecm
\centerline{\bf On the Integral Cohomology of Wreath Products}
\bigskip
\centerline{Ian J. Leary}
\bigskip
\centerline{\it Faculty of Mathematical Studies, University of Southampton}
\centerline{\it Southampton, SO17 1BJ}
{{\narrower\smallskip\noindent
\bigskip
{\smallbold Abstract.} {\smallfont Under mild conditions on the space
$\scriptstyle X$, we describe the additive structure
of the integral cohomology of the space 
$\scriptstyle X^p\times_{C_p}EC_p$ in terms of the cohomology of 
$\scriptstyle X$.  
We give weaker results for other similar spaces, and deduce various
corollaries concerning the cohomology of finite groups.}
\smallskip}}
\bigskip

\beginsection 0. INTRODUCTION

Let $S$ be a group with a fixed action on a finite set $\Omega$.  By
the wreath product $G\wr S$ of a group $G$ with $S$ we mean a split
extension with kernel $G^\Omega$, quotient $S$, and with the
$S$-action on $G^\Omega$ given by permuting the copies of $G$.  Our
main interest is the integral cohomology of finite groups of the form
$G\wr S$.  We work in greater generality however, because it is no
more difficult to study the cohomology of spaces of the form
$X^\Omega\times_SE$, where $E$ is an $S$-free $S$-CW-complex, and $X$
is a CW-complex of finite type.  The mod-$p$ cohomology of certain
such spaces plays a crucial role in Steenrod's definition of the 
reduced power operations [\ste].  Building on work of Steenrod,
Nakaoka described the cohomology of such spaces with coefficients in
any field [\nak]. 
The point about working over a field is that then the cellular cochain
complex for $X$ is homotopy equivalent to the cohomology of $X$,
viewed as a complex with trivial differential.  If the integral
cohomology of $X$ is free, then a similar result holds in this case.
Evens used this to study the cohomology of (the classifying space of) 
the Lie group $U(m) \wr \Sigma_n$ in the course of his
work on Chern classes of induced representations [\eva].  
The study of the integral cohomology in the case when $H^*(X)$ is not
free is much harder.  The pioneers in this case were Evens and Kahn,
who made a partial study of the important special case of
$X^p\times_{C_p}EC_p$.  We complete the study of this case in
Section~4 below, which could be viewed as both an extension of and a
simplification of [\evk, section~4].  Many, but not all, of our results are
corollaries of this work.  Our paper has the following structure.  
 
%\headline{\sevenrm\hfil COHOMOLOGY OF WREATH PRODUCTS\hfil}
In Section~1 we give some algebraic background.  Most of this material
is well-known, although we have not seen Lemma~1.4 stated explicitly
before, and we believe that Lemma~1.1 is original.  This lemma, which
compares spectral sequences coming from double complexes consisting of
\lq the same groups', but with \lq different differentials', is the
key to our extension of Evens-Kahn's results [\evk].  Theorem~2.1 is a
statement of the result of Nakaoka mentioned above, which for
interest's sake we have made more general than the original.  A weak
version of this theorem could be stated as: \lq if $H^*(X)$ is free,
then the Cartan-Leray spectral sequence for $H^*(X^\Omega\times_SE)$
collapses at the $E_2$-page'.  A similarly weakened version of
Theorem~2.2 would say that over the integers (or any PID), this
spectral sequence collapses at the $E_r$-page, $r= 2 + |\Omega|$,
without any condition on $H^*(X)$.  

The Cartan-Leray spectral sequence may be obtained from a  double
complex.  In Sections 3~and~4 we study the other spectral sequence
associated with the same double complex.  In Section~3 we define the
elements $\alpha\wr 1$ and give upper and lower bounds on their
orders.  In Section~4 we describe the additive structure of the
integral cohomology of $X^p\times_{C_p}EC_p$, extending work of
Evens-Kahn [\evk].  One surprising corollary of this result is that
for $p\geq 5$, a cyclic summand of $H^i(X)$ of order $p^j$ gives rise
to $(p-1)/2$ cyclic summands of $H^*(X^p\times_{C_p}EC_p)$ of order
$p^{j+1}$, not just the summand that one would expect in degree $pi$.
In Section~5 we use the results of Section~4 to determine the
differentials and extension problems in the Cartan-Leray spectral
sequence for $X^p\times_{C_p}EC_p$.  

Quillen's detection lemma [\qui], which states that the mod-$p$
cohomology of $X^p\times_{C_p}EC_p$ is detected by two maps 
from $X^p$ and $X\times BC_p$, is a corollary of Nakaoka's results.
In Section~6 we describe the kernel of the analogous map in integral
cohomology.  In Section~7 we sketch how the methods of Sections~4,
5,~and~6 may be applied to describe the $p$-local cohomology of
$X^p\times_{\Sigma_p}E\Sigma_p$.  In fact this turns out to be
considerably simpler than the case of of the cyclic group.  

In Section~8 we review a conjecture of A.~Adem and H.-W.~Henn
concerning the exponent of integral cohomology of finite groups.  We
show that if $G$ does not afford a counterexample, then neither does
$G\wr C_p$.  For a slightly stronger conjecture we obtain the more
general result that if $G$ and $S$ do not afford counterexamples then
neither does $G\wr S$.  We also recall an upper bound on the
(eventual) exponent of integral cohomology of $p$-groups, based on
generalized Frattini subgroups, that we gave in [\ian].  Section~9 
describes an example showing that this bound is not always best
possible.  The only connection between this section and wreath
products is that evidence gleaned from wreath products led the author
to believe for some time that such examples could not exist.   

\headline{\sevenrm\hfil COHOMOLOGY OF WREATH PRODUCTS\hfil}
In Section~10, we assume that the $p$-local cohomology of $X$ is a
finitely generated algebra, so that the variety of all ring
homomorphisms from $H^*(X)$ to an algebraically closed field $k$ of
characteristic~$p$ is affine.  In this case the variety for
$H^*(X^p\times_{C_p}EC_p)$ is also affine, and may be described in
terms of the variety for $H^*(X)$.  (This is easily deduced from
[\nak], but is first stated explicitly in [\qui].)  We use the results
of Sections 4~and~6 to give a similar description, for each $i$, of
the subvariety corresponding to the annihilator in
$H^*(X^p\times_{C_p}EC_p)$ of the element $p^i$.  First we prove some
general properties of such subvarieties.  In Section~11 we apply the
results of Section~10 to cohomology of finite groups.  The
subvarieties that we study were introduced in this context by Carlson
[\carl].  Let $G$ be a $p$-group, let $W(G)$ be the variety of ring
homomorphisms from $H^*(G)$ to $k$, and let $W_i(G)$ be the subvariety
corresponding to $\Ann(p^i)$.  Each $W_i(G)$ is a covariant functor of
$G$.  We show that $W_i(G)$ is contained in the image of
$W(\Phi_i(G))$, where $\Phi_i(G)$ is the generalized Frattini subgroup
introduced in Section~8.  Carlson asked in [\carl] if the image of
$W(\Phi^i(G))$ is contained in $W_i(G)$, where $\Phi^i(G)$ is the
$i$th iterated Frattini subgroup of $G$.  We give examples where this
does not hold.  We show however that the image of $W(\Phi(G)\cap
Z(G))$ is contained in $W_1(G)$.  We close with the remark that in all
known examples, each $W_i(G)$ has a nice description.  

\noindent
{\bf Acknowledgements.}  The author worked part-time on this
material during a period in which he held post-doctoral fellowships at
the ETH, Z\"urich, the CRM, Barcelona, and the MPI, Bonn, funded by
the ETH, the DGICYT, and the Leibniz fellowship scheme respectively.  
The author gratefully acknowledges the support and hospitality of
these institutions.  The author thanks Guido Mislin and Urs Stammbach 
for their advice and encouragement during the early stages of this
work, and Jon Carlson and Leonard Evens for helpful email
correspondence.

\beginsection 1. NOTATION AND ALGEBRAIC PRELIMINARIES

We hope that our notation is either defined where used or is evident.
As a rough guide we mention the following points.  

\mkpt
$p$ is a prime number throughout.  

\mkpt
$C_n$ is a cyclic group of order $n$.  

\mkpt
$\Omega$ is the set $\{1,\ldots,l\}$, or occasionally
$\{1,\ldots,p\}$.  

\mkpt
$\Sigma_l=\Sigma(\Omega)$, is the symmetric group on $\Omega$.  

\mkpt
All spaces are taken to be CW-complexes, and the topology on a product
is chosen so that the product of CW-complexes is a CW-complex in the
natural way.   

\mkpt
For $G$ a group, $EG$ is a contractible $G$-free $G$-CW-complex.  

\mkpt
$R$ is a commutative ring, often the integers $\Bbb Z$ or the
$p$-local integers $\Bbb Z_{(p)}$.  Usually $-\otimes-$ and
$\Hom(-,-)$ should be taken over $R$.  In Sections 10~and~11 however,
$R$ is a more general commutative $\plocal$-algebra, and $-\otimes-$
stands for the tensor product over $\plocal$.  

\mkpt 
Usually $S$ stands for a subgroup of $\Sigma(\Omega)$, although in
Section~9 a sphere called $S$ makes a brief appearance, and in
Section~10 $S$ is a commutative $\plocal$-algebra.  

\mkpt
All chain complexes and cochain complexes are bounded below.  

Double cochain complexes will be denoted $(E_0^{*,*},d,d')$,
$E_0^{*,*}$, or similarly, where the subscript is intended 
to suggest that we will
take spectral sequences.  Double complexes (and pages of the
associated spectral sequences) will be depicted with the first index
running horizontally: 
$$\matrix{
&\vdots&&\vdots&\cr 
\cdots&E_0^{i,j+1}&\mapright{d}&E_0^{i+1,j+1}&\cdots\cr
&\mapup{d'}&&\mapup{d'}&\cr
\cdots&E_0^{i,j}&\mapright{d}&E_0^{i+1,j}&\cdots\cr
&\vdots&&\vdots&\cr }.$$
Following Cartan-Eilenberg [\caei], we call the associated spectral
sequence in which $d_0=d'$ the \lq type~I spectral sequence', and the
spectral sequence in which $d_0=d$ the \lq type~II spectral sequence'.
Note that higher differentials in type~II spectral sequences will
point upwards and leftwards on our illustrations.  

The following lemma seems to be new, and is very useful in our
calculation (in Section~4) of $H^*(X^p\times_{C_p}EC_p)$.

\proclaim Lemma 1.1.  Let $(E_0^{*,*},d,d')$ be a double cochain
complex of abelian groups, and let $E_r^{*,*}$ be the
corresponding type II spectral sequence.  
For any integer $n$, make a second 
double cochain complex $(\tilde E_0^{*,*},\tilde d, \tilde d')$,
where $\tilde E_0^{i,j}=E_0^{i,j}$, $\tilde d=d$, and $\tilde d'
= nd'$.  Then any element $x$ of $E_0^{*,*}$ that represents an
element of $E_r^{*,*}$ also represents an element of 
$\tilde E_r^{*,*}$, and $\tilde d_r(x)= n^rd_r(x)$.  

\proof By assumption, there exist $x=x_1,\ldots,x_r\in E_0^{*,*}$
such that $d(x_1)=0$ and $d(x_{i+1})+d'(x_i)=0$ for $1\leq i< r$,
and then by definition $d_r(x)$ is represented by $d'(x_r)$.  
Now if we let $y_i= n^{i-1}x_i$, we see that $x=y_1$, 
$\tilde d(y_1)=0$, $\tilde d(y_{i+1})+\tilde d'(y_i)=0$ for 
$1\leq i<r$, and that 
$\tilde d_r(x)=\tilde d'(y_r)=n^rd'(x_r)=n^rd_r(x)$. 
\qed

\remark There is no simple relation between the type I spectral
sequences associated with the complexes $E_0^{*,*}$ and $\tilde
E_0^{*,*}$.  

\proclaim Proposition 1.2. If $C^*$ is a cochain complex of
$R$-modules, then $C^{*\otimes\Omega}$ may be given the structure of a
cochain complex of $R\Sigma(\Omega)$-modules as follows.  The action
of the transposition $\sigma_i=(i,i+1)\in \Sigma(\Omega)$ is defined
on a homogeneous element $c_1\otimes\cdots\otimes c_l$ by 
$$\sigma_i(c_1\otimes\cdots\otimes c_l)= 
(-1)^{\deg(c_i)\deg(c_{i+1})}c_1\otimes\cdots\otimes 
c_{i+1}\otimes c_i\otimes\cdots\otimes c_l.$$

\proof First check that the action of each $\sigma_i$ commutes with
the differential on $C^{*\otimes\Omega}$.  Now recall that
$\Sigma(\Omega)$ has the following presentation as a Coxeter group
[\brobui]:  
$$\Sigma(\Omega)= \langle \sigma_1,\ldots,\sigma_{l-1} \mid 
\sigma_i^2, (\sigma_i\sigma_j)^2, (\sigma_k\sigma_{k+1})^3\rangle, $$
where $1\leq i,j < l$, $1\leq k < l-1$, and $|i-j|>1$.  Now check that
each of the relators in this presentation acts trivially on
$C^{*\otimes \Omega}$.  \qed

\remark The above proof was suggested to the author by Warren Dicks.
It seems to be easier than proofs in which one works out the sign for
the action of every element of $\Sigma(\Omega)$  (see for example 
[\ben,\evb]).  On the other hand, one needs to know the sign of 
the action of an arbitrary element in most applications.  

If $C^*$ is the cellular cochain complex on a finite type CW-complex
$X$, it may be shown that the cellular action of $\Sigma(\Omega)$ on
$X^\Omega$ (with the product cell structure) induces the above action
on $C^*(X^\Omega)= C^*(X)^{\otimes \Omega}$.  

Let $S$ be a subgroup of $\Sigma(\Omega)$, and let $W_*$ be a chain
complex of free $RS$-modules.  
Much of our work consists of studying double cochain complexes of the
form 
$$E_0^{i,j}= \Hom_S(W_i, (C^{*\otimes\Omega})^j),$$
and their associated spectral sequences.  In many cases, $W_*$ will
be either the cellular cochain complex of a free $S$-CW-complex $E$,
or will be a free resolution of the trivial $RS$-module $R$.  In
either of these cases there is a homomorphism (of $RS$-complexes,
where $S$ has the diagonal action on $W_*\otimes W_*$):  
$$\Delta: W_*\mapright{} W_*\otimes W_*$$
which gives rise to an anticommutative ring structure on
$H^*\Hom_S(W_*,R)$.  In this case, the isomorphism 
$(C^{*\otimes\Omega})\otimes R\cong (C^{*\otimes\Omega})$ gives rise
to a graded $H^*\Hom_S(W_*,R)$-module structure on $H^*\tot
E_0^{*,*}$, and a bigraded action on each of the two spectral 
sequences from the $E_2$-page onwards.  (In each case, 
$H^*\Hom_S(W_*,R)$ is graded in the $i$-direction.)  

Lemma~1.3, which we shall state but not prove, is due to Steenrod
[\ste].  It is used in Steenrod's definition of the reduced
$p$-powers, as well as in Nakaoka's work on wreath products (see 
[\nak], or Theorem~3.1 below).  For most of our purposes the easier
Lemma~1.4 will suffice however.  Lemma~1.4 is implicit in a remark in
[\eva], but I know of no explicit statement of it in the literature.

\proclaim Lemma 1.3.  If $W_*$ is (as above) a complex of free
$RS$-modules and $f,g:C'^*\rightarrow C^*$ are homotopic maps between
cochain complexes of $R$-modules, then the maps $\Hom_S(1,f^{\otimes\Omega})$
and $\Hom_S(1,g^{\otimes\Omega})$ from
$\Hom_S(W_*,C'^{*\otimes\Omega})$ to $\Hom_S(W_*,C^{*\otimes\Omega})$ 
are homotopic.  \qed

The proof of Lemma 1.3 is slightly simpler if it is also assumed that
$W_*$ is acyclic, and this is the statement given in the recent books
[\ben] and [\evb].  Even if one is only interested in the special
case of group cohomology however, it is worth having the stronger
version of Lemma~1.3 because this justifies the application of
Nakaoka's argument to wreath products of the form $G\wreath S'$,
where $S'$ has $S$ as a quotient and acts on $G^\Omega$ via the
action of $S$ on $\Omega$.  

\proclaim Lemma 1.4.  Let $W_*$ be a complex of free $RS$-modules, and
let $f:C'^*\rightarrow C^*$ be a homotopy equivalence between cochain
complexes of $R$-modules.  Define a double cochain complex $E_0^{*,*}$
by 
$$E_0^{i,j}=\Hom_S(W_i,(C^{*\otimes\Omega})^j),$$
and similarly for $E_0'^{i,j}$.  Then the map 
$$\Hom_S(1,f^{\otimes\Omega}):E_0'^{*,*}\rightarrow E_0^{*,*}$$
of double cochain complexes induces an isomorphism on the $E_1$-pages
of the associated type~I spectral sequences, and hence an isomorphism
between the homologies of the corresponding total complexes.  

\proof The map $f^{\otimes\Omega}:C'^{\otimes\Omega}\rightarrow 
C^{\otimes\Omega}$ is a homotopy equivalence of $R$-complexes, and is
an $RS$-map, so induces an $RS$-module isomorphism between 
$H^j(C'^{\otimes\Omega})$ and $H^j(C^{\otimes\Omega})$.  If $W_i$ is
free with basis $I_i$, then $E_1^{i,j}$ (resp.\ $E_1'^{i,j}$) is
isomorphic to a product of copies of $H^j(C^{\otimes\Omega})$ 
(resp.\ $H^j(C'^{\otimes\Omega})$) indexed by $I_i$, and the map is
an isomorphism as claimed.  \qed

The conclusion of Lemma 1.4 does not necessarily hold for $f:C'^*
\rightarrow C^*$ that only induces an isomorphism on cohomology, but 
recall the following lemma, from for example [\spa, p167].  

\proclaim Lemma 1.5.  A map between $R$-projective cochain complexes is
a homotopy equivalence if and only if it induces an isomorphism on
cohomology.  \qed

We end this section with some remarks concerning the case when $R$ is
a principal ideal domain (PID) and $C^*$ is an $R$-free cochain
complex such that each $H^i(C^*)$ is a finitely generated $R$-module.
For each $i$, fix a splitting of $H^i(C^*)$ as a direct sum of cyclic
$R$-modules.  Thus for some indexing set $A$, 
$$H^*(C)\cong \bigoplus_{a\in A}H(a),$$
where $H(a)$ is isomorphic to $R/(r(a))$ and is a summand in degree
$i(a)$, for some $r(a)\in R$ and integer $i(a)$.  

For $r\in R$ and $i\in \Bbb Z$, define a cochain complex $C(r,i)^*$
such that each $C(r,i)^j=0$, except that $C(r,i)^i\cong R$ and if 
$r$ is non-zero, then $C(r,i)^{i-1}\cong R$.  Define the differential
in $C(r,i)^*$ so that $H^i(C(r,i)^*)\cong R/(r)$.  Define a cochain
complex $C'^*$ by 
$$C'^*=\bigoplus_{a\in A} C(a)^*,$$
where $C(a)^*=C(r(a),i(a))$.  It is easy to construct a map from
$C'^*$ to $C^*$ inducing an isomorphism on cohomology, which is
therefore a homotopy equivalence by Lemma~1.5.  

Now consider the $RS$-module structure of the complex
$(C'^{*\otimes\Omega})^*$.  (The motivation for this is provided by Lemmata 
1.3~and~1.4.)  As a complex of $R$-modules, 
$(C'^{*\otimes\Omega})^*$ splits as a direct sum of pieces of the form 
$$C(a_1,\ldots,a_l)= C(a_1)^*\otimes\cdots\otimes C(a_l)^*.$$
Call such a summand a \lq cube', and call the $C(a_j)^*$ that arise
the cube's \lq sides'.  The action of $S$ permutes the cubes, and the
stabilizer $S'$ of the cube $C(a_1,\ldots,a_l)^*$ is equal to the
stabilizer of $(a_1,\ldots,a_l)\in A^\Omega$.  Each $R$-module summand of 
$C(a_1,\ldots,a_l)^*$ of the form 
$C(a_1)^{i(a_1)-\epsilon(1)}
\otimes\cdots\otimes C(a_l)^{i(a_l)-\epsilon(l)}$, 
where $\epsilon(j)=
0$~or~1 and $\epsilon(j)=0$ if $r(a_j)=0$, is in fact an
$RS'$-summand.  Define $\Omega^-$ by 
$$\Omega^-= \{j\in \Omega\mid i(a_j)-\epsilon(j) \hbox{ is odd }\},$$
and $\Omega^+= \Omega\setminus \Omega^-$.  Then $S'$ is a subgroup of 
$\Sigma(\Omega^+)\times \Sigma(\Omega^-)\leq \Sigma(\Omega)$, and 
the action of $S'$ on $C(a_1)^{i(a_1)-\epsilon(1)}
\otimes\cdots\otimes C(a_l)^{i(a_l)-\epsilon(l)}$ is given by the sign
action of $\Sigma(\Omega^-)$ tensored with the trivial action of 
$\Sigma(\Omega^+)$.  (This may be checked using Proposition~1.2.)

\beginsection 2. NAKAOKA'S ARGUMENT AND GENERALIZATIONS

Recall that $S$ is a subgroup of the symmetric group 
on a finite set $\Omega$, and 
let $E$ be a free $S$-CW-complex of finite type.  Note that 
we do not require $E$ to be contractible.  For any finite type CW-complex 
$X$ we may define an action of $S$ on a product of copies of $X$ indexed 
by $\Omega$, letting $S$ permute the factors.  Now $S$ acts \lq 
diagonally' on $X^\Omega\times E$, and we may form the quotient space 
$(X^\Omega\times E)/S$.  Topologists usually insist that $S$ should act 
on the right of $X^\Omega$ and on the left of $E$, and write $X^\Omega 
\times_SE$ for the quotient.  Assume that $E$ is connected.  
By considering various covering spaces of 
$X^\Omega\times_SE$ it is easy to see that in this case 
the fundamental group of $X^\Omega\times_SE$ is the  
wreath product $G\wreath \pi_1(E/S)$, where $G$ is the fundamental group
of $X$.  By this wreath product we mean a split extension with kernel 
$G^\Omega$ and quotient $\pi_1(E/S)$, where $\pi_1(E/S)$ acts via the 
action of its quotient $S$ on $\Omega$.  In particular, if $E$ is simply 
connected, then $\pi_1(X^\Omega\times_SE)$ is the wreath product 
$G\wreath S$.  From now on we shall assume that $E$ is both connected
and simply connected for clarity---it is easy to make the necessary 
changes to cover the general case.  

For any field $k$, Nakaoka's argument determines $H^*(X^\Omega\times_SE;k)$
in terms of $H^*(X;k)$ and $H^*(E/S;N)$ for various signed permutation 
$kS$-modules $N$.  The argument applies equally to the cohomology of 
\def\symsp{{X^\Omega\times_SE}}
$\symsp$ with non-trivial coefficients coming from some 
$kG\wreath S$-modules (in fact, those which are \lq tensor-induced' 
from $kG$-modules) and gives some information about $H^*(\symsp;R)$ for 
any commutative ring $R$.  The following account includes both of these 
generalizations.  

Let $R$ be a commutative ring, and let $\Omega$, $S$, $E$ and $X$ be as 
above, with $\pi_1(X)=G$.  All chain complexes will be chain complexes 
of $R$-modules unless otherwise stated and all unmarked tensor products 
will be over $R$.  Let $W_*$ be the cellular $R$-chain complex for $E$, 
so that $W_*$ is a chain complex of free $RS$-modules, and is a free 
resolution for $R$ over $RS$ in the case when $E$ is acyclic.  Let $U_*$ 
be the cellular chain complex of the universal cover of $X$, so that 
$U_*$ is a chain complex of finitely generated free $RG$-modules.  If 
$M$ is an $RG$-module, then $RG^\Omega$ acts on $M^{\otimes\Omega}$ by 
letting $G^\Omega$ act component-wise, and $S$ acts on $M^{\otimes\Omega}$
by permuting the factors of each monomial $m_1\otimes m_2\otimes\ldots
\otimes m_l$.  
Together these actions define an $RG\wreath S$-module structure on 
$M^{\otimes\Omega}$.  Similarly, the total complex of $U_*^{\otimes\Omega}$
becomes a complex of $RG\wreath S$-modules, where the action of $S$ is
as described in Proposition~1.2.  From this point of view the signs
that arise are due to the fact that the action of $S$ on the 
$\Omega$-cube does not necessarily preserve its orientation.  The
complex $U_*^{\otimes\Omega}$ is the cellular chain complex of the
universal cover of $X^\Omega$ (with the product CW structure), where
the action of $G\wreath S$ comes from the action of $S$ on $X^\Omega$.
Similarly, the complex $U_*^{\otimes\Omega}\otimes W_*$, with the
diagonal action of $G\wreath S$, is eaily seen to be the cellular
chain complex of the universal cover of $\symsp$.  

Define a double cochain complex $E_0^{i,j}$ by: 
$$ E_0^{i,j}=\Hom_{RG\wreath S}((U_*^{\otimes\Omega})_j\otimes W_i,
M^{\otimes\Omega}).$$
The associated total complex is the complex of $G\wreath
S$-equivariant cochains on the universal cover of $\symsp$ with
values in $M^{\otimes\Omega}$, and so the cohomology of this complex
is $H^*(\symsp;M^{\otimes\Omega})$.  There are two spectral sequences
associated to the double complex $E_0^{i,j}$.  The type~I spectral
sequence has $E_2$-page as follows:  
$$E_2^{i,j}=H^i(E/S;H^j(X^\Omega;M^{\otimes\Omega})).$$
This is a spectral sequence of Cartan-Leray type for the covering 
of $\symsp$ by $X^\Omega\times E$.  From the $E_2$-page 
onwards it is isomorphic to the
Serre spectral sequence for the fibration 
$$X^\Omega \longrightarrow \symsp \longrightarrow E/S.$$
Each of the two spectral sequences admits a 
bigraded $H^*(E/S;R)$-module structure from the $E_2$-page onwards, as
was shown in Section~1.  

So far we have not  used the assumption that $X$ is of finite type.  
This is required to establish the second of the following
isomorphisms of double cochain complexes:  
$$\eqalign{E_0^{*,*}=&\Hom_{G\wreath S}(U_*^{\otimes\Omega}\otimes W_*,
M^{\otimes\Omega})\cr
\cong&\Hom_S(W_*,\Hom_{G^\Omega}(U_*^{\otimes\Omega},M^{\otimes\Omega}))\cr
\cong&\Hom_S(W_*,\Hom_G(U_*,M)^{\otimes\Omega}).\cr}$$
Thus we are reduced to the study of the cohomology of a double complex
of the type discussed in Section~1.  

\remark
(The algebraic case.)  The case of interest in group cohomology may be
recovered as the case when $X$ and $E/S$ both have trivial higher
homotopy groups (note that this is
more general than requiring $E$ to be contractible).  Of course, in 
the algebraic case we can do without the space $X$ altogether, and
just take the complex $U_*$ to be a finite type ($R$-free) $RG$-projective
resolution for $R$.  Thus the above argument applies 
to groups $G$ of type $FP(\infty)$ over $R$, rather than just those
groups $G$ having a $K(G,1)$ of finite type.  

The following theorem is due to Nakaoka, although we have deliberately
made the statement more general than the original (which considered
only trivial coefficients).  
  
\proclaim Theorem 2.1 ([\nak,\eva,\kah]).  As above, let $X$ be a
connected CW-complex of finite type with fundamental group $G$, and
let $M$ be an $R$-free $RG$-module such that $H^*(X;M)$ is
$R$-projective.  Let $E$ be a (connected) free $S$-CW-complex.  Then
there is an isomorphism of graded $R$-modules as follows.  
$$H^*(X^{\otimes\Omega}\times_SE;M^{\otimes\Omega})\cong
H^*(E/S;H^*(X;M)^{\otimes\Omega})$$
The spectral sequence with coefficients in $M^{\otimes\Omega}$ 
for the fibration $X^\Omega\rightarrow X^\Omega\times_SE\rightarrow 
E/S$ collapses at the $E_2$-page, and the additive extensions in the 
reconstruction of the cohomology from the $E_\infty$-page are all
trivial.  

\remarks
Note that the freeness conditions are automatically satisfied if $R$
is a field.  In the case when $M=R$, the isomorphism of the theorem
is moreover an isomorphism of $R$-algebras, and for general $M$ it
may be shown to be an isomorphism of $H^*(X^\Omega\times_SE;R) \cong
H^*(E/S;H^*(X;R)^{\otimes\Omega})$-modules.  

\proof Under the conditions of the theorem, both $C^*=\Hom_G(U_*,M)$
and $C'^*=H^*(C^*)$ are $R$-projective.  Viewing $C'^*$ as a cochain
complex with trivial differential it is easy to construct a chain map
from $C'^*$ to $C^*$ inducing the identity map on cohomology (see [\nak]
or [\evk]), which is a homotopy equivalence by Lemma~1.5.  Now by
either Lemma~1.3 or Lemma~1.4 we see that
$H^*(X^\Omega\times_SE;M^{\otimes\Omega})$ may be calculated as the
cohomology of the total complex of
$E_0'^{i,j}=\Hom_S(W_i,(C'^{*\otimes\Omega})^j)$.  But this double
complex has trivial $j$-differential, so splits as a direct sum of
double complexes concentrated in constant $j$-degree (or \lq rows').
\qed

\remarks 
There are other proofs of special cases of Theorem~2.1.  There is a
topological proof due to Adem-Milgram in the special case when $S$ 
is cyclic of order~$p$, $E$ is contractible, $S$ acts freely 
transitively on $\Omega$, and the coefficients form a field [\ademil].  

There is also a short proof of the trivial coefficient case arising 
in cohomology of finite groups, due to Benson and Evens ([\ben,
vol.~II, p.130 and \evb, Thrm.~5.3.1]).  
This proof is as follows:  If $G$ is finite, and $R=k$ is 
a field, then one may take $U_*$ to be 
a minimal resolution for $k$ over $kG$.  In this case, if $M$ is a simple
$kG$-module, for example the trivial module $k$, then the differential
in $\Hom_G(U_*,M)$ is trivial
and so the double complex for computing $H^*(G\wreath S;
M^{\otimes\Omega})$ has one of its differentials trivial, and the
type~I spectral sequence collapses at the $E_2$-page.  
It is not true however that if also $W_*$ is a minimal resolution for
$k$ over $kS$ then $U_*^{\otimes\Omega}\otimes W_*$ is a minimal
resolution for $k$ over $kG\wreath S$, because in general this
will have a larger growth rate than a minimal resolution.  

If the condition on $H^*(X;M)$ is weakened, Nakaoka's argument may
still be applied, but the conclusion is far weaker.  For example,
consider the following theorem.  

\proclaim Theorem 2.2.  Take notation and conditions as in the
statement of Theorem 2.1, but replace the condition that $H^*(X;M)$
should be $R$-projective by the condition that $H^*(X;M)$ should have
projective dimension at most one over $R$.  Then from $E_2$ onwards, 
the spectral sequence for the fibration 
$X^{\Omega}\rightarrow X^\Omega\times_SE\rightarrow E/S$
with $M^{\otimes\Omega}$ coefficients is a direct sum of spectral
sequences $E_r^{*,*}=\bigoplus_\alpha E_{r,\alpha}^{*,*}$, where 
each $E_{r,\alpha}^{i,j}$ has \lq height' at most $|\Omega|$.  There
is a corresponding direct sum decomposition of
$H^*(X^\Omega\times_SE;M^{\otimes\Omega})$.  

\remark By the phrase \lq $E_{r,\alpha}^{*,*}$ has height at most
$|\Omega|$' we mean that there exist integers $n(\alpha)\leq
n'(\alpha)$ such that $E_{r,\alpha}^{i,j}=0$ if $j<n(\alpha)$ or 
$j>n'(\alpha)$, and $n'(\alpha)-n(\alpha)\leq |\Omega|$.  In
particular this implies that $d_{|\Omega|+1}$ is the last possibly
non-zero differential.  

\proof Let $P_i\inj Q_i\surj
H^i(X;M)$ be a projective resolution over $R$ for
$H^i(X;M)$, and build a cochain complex $C'^*$ with $C'^i=Q_i\oplus
P_{i+1}$ and differential given by the following composite.  
$$C'^i\surj P_{i+1}\inj Q_{i+1}\inj C'^{i+1}$$
Then using the projectivity of $P_i$ and $Q_i$ it is easy to
construct a chain map from $C'^*$ to $\Hom_G(U_*,M)$ inducing an
isomorphism on cohomology, which is an equivalence by Lemma 1.5.  Now
by either Lemma 1.3 or Lemma 1.4, the double complex $E_0'^{*,*}= 
\Hom_S(W_*,C'^{\otimes\Omega})$ may be used to compute
$H^*(X^\Omega\times_SE;R)$.  By construction $C'^*$ splits as a
direct sum of complexes of length at most one, and so
$(C'^{\otimes\Omega})^*$ splits as a direct sum of complexes of
length at most $|\Omega|$, and hence $E_0'^{*,*}$ splits as a direct
sum of double complexes of height at most $|\Omega|$.  The claimed
properties of the spectral sequence now follow from Lemma 1.4.  
\qed

\remarks There is no easy generalization of the statement given after
Theorem~2.1 concerning the ring structure of
$H^*(X^\Omega\times_SE;R)$.  It is known, for example, that the ring
structure of the integral cohomology of $X$ does not suffice to
determine that of $X\times X$ [\pal].  This could be considered as
the case of $X^\Omega\times_SE$ when $R$ is the integers, $E$ is a
point, and $S$ is the trivial subgroup of the symmetric group on a
set $\Omega$ of size two.  

One may ask about similar results to the above for generalized
cohomology theories.  We believe that the following statement is
implicit in [\lmm, chap.~IX]: \lq Let $h^*(-)$ be a complex oriented
generalized cohomology theory
associated with an ${\rm H}_\infty$-ring spectrum, and assume that $h^*(X)$
is free over $h^*$ with basis $B$ concentrated in even degrees.  Then
$h^*(X^\Omega\times_SE)$ is isomorphic to a direct sum of 
$\tilde h^*(E/S')$'s, shifted in degree, where $S'$ runs over the stabilizers
of a set of orbit representatives for the action of $S$ on
$B^\Omega$.'  For specific choices of $S$ and $E$ stronger results are
known---see [\jrh] for some recent results in the case when $S=C_p$ and $E$
is contractible.  

\beginsection 3.  THE ORDER OF $\wre \alpha$  

Throughout this section, take $R$ to be either $\Bbb Z$ or one of its
localizations, and let $W_*$ be the cellular chain complex of an
$S$-free $S$-CW-complex $E$ with augmentation $\epsilon:
W_*\rightarrow R$, so that the homology of $\ker(\epsilon)$ is the
reduced $R$-homology of $E$.  For any $X$ of finite type and any
$\alpha \in H^{2i}(X;R)$, $\wre \alpha\in
H^{2il}(X^\Omega\times_SE;R)$ may be defined
([\ste,\nak,\evb]).  We give bounds on the
order of $\wre \alpha$ in terms of the order of $\alpha$.  In fact it
costs no extra work to replace $C^*(X;R)$ by an arbitrary $R$-free
cochain complex of finite type, so we do so.  Much of the section
generalizes to the case when $R$ is any PID, if the statements
concerning orders are replaced by statements about annihilators.  

For $C^*$ a finite type $R$-free cochain complex, and $c\in C^*$ of
even degree, define $\wre c \in \Hom_S(W_*,C^{*\otimes\Omega})$ by 
$$\wre c (w) = \epsilon(w).c\otimes\cdots\otimes c, $$
for any $w\in W_*$.  The same formula may be used for $c$ of odd
degree, but in this case it defines an element of
$\Hom_S(W_*,C^{*\otimes\Omega}\otimes \widehat R)$, where $\widehat R$
is the $RS$-module of $R$-rank one on which $S$ acts via the sign
representation of $\Sigma(\Omega)$.  If $c$ is a cocycle, then so is
$\wre c$, and using Lemma~1.3 it may be shown that the cohomology
class of $\wre c$ depends only on the cohomology class of $c$ (this is
the only place where we require Lemma~1.3 rather than Lemma~1.4).  
Thus if $\alpha \in H^*(C^*)$, we may define a unique element 
$\wre \alpha$ of $H^*\tot\Hom_S(W_*,C^{*\otimes\Omega})$.  

If $\alpha$ does not generate a direct summand of $H^*(C^*)$, pick
$\alpha'$ and an integer $n$ such that $\alpha'$ generates a direct
summand of $H^*(C^*)$ and $\alpha=n\alpha'$.  Then (check on the level
of cochains) $\wre\alpha = \wre{(n\alpha')} = n^l(\wre \alpha)$.
Hence it suffices to consider the case when $\alpha$ generates a 
summand of $H^*(C^*)$.  

\proclaim Theorem 3.1.  Fix $S$, $W_*$, $R$, as above.  Let $C^*$ be
any $R$-free finite type cochain complex, and let $\alpha \in
H^*(C^*)$.  Then the order of $\wre \alpha\in
H^*\tot\Hom_S(W_*,C^{*\otimes\Omega})$ depends only on the order
$O(\alpha)$ of $\alpha$, and the order of a cyclic summand of
$H^*(C^*)$ containing $\alpha$.  If $\alpha$ has infinite order, then
so does $\wre \alpha$.  If $\alpha$ generates a finite direct summand
of $H^*(C^*)$, then 
$$O(\alpha)\leq O(\wre\alpha) \leq l'.O(\alpha),$$
where $l'$ is the h.c.f.\ of the lengths of the $S$-orbits in
$\Omega$.  

\proof The reduction to the case when $\alpha$ generates a direct
summand of $H^*(C^*)$ is discussed above.  If $\alpha$ does generate a
direct summand of $H^*(C^*)$, pick a splitting of $H^*(C^*)$ into
cyclic summands as at the end of Section~1:  
$$H^*(C^*)=\bigoplus_{a\in A}H(a).$$
Choose the splitting in such a way that $\alpha$ generates one of the
summands, say $H(a_0)$.  Construct $C'^*$ as in Section~1.  Then 
$H^*(C(a_0)^*)$ maps isomorphically to the summand of $H^*(C^*)$
generated by $\alpha$, and $C(a_0)^{*\otimes\Omega}$ is an $RS$-summand of
the complex $C'^{*\otimes\Omega}$. Hence the image of
$H^*\tot\Hom_S(W_*,C(a_0)^{*\otimes\Omega})$ in
$H^*\tot\Hom_S(W_*,C^{*\otimes\Omega})$ is a direct summand.  Moreover, if
$\alpha'\in H^*(C(a_0)^*)$ maps to $\alpha$, then $\wre{\alpha'}$ maps
to $\wre\alpha$.  Thus we see that if $\alpha$ generates a summand of
$H^i(C^*)$ of order $n$, then $O(\wre\alpha)= O(\wre{\alpha'})$, where
$\alpha'$ generates $H^i(C(n,i)^*)$ (see Section~1 for the definition of
$C(n,i)^*$).  Let $D(n,i)^*= C(n,i)^{*\otimes\Omega}$.  It may be
checked that as complexes of $RS$-modules, for any $i$, 
$$D(n,2i)^{*-2li}\cong D(n,0)^*\cong D(n,2i+1)^{*-2li-l}\otimes \widehat R,$$
where $\widehat R$ is the sign representation of $S\leq
\Sigma(\Omega)$.  Hence the order of $\wre {\alpha'}$ does not depend
on $i$.  

Let $E'$ be a set of points permuted freely, transitively by $S$, and
let $E''$ be a contractible free $S$-CW-complex with one orbit of zero
cells.  Then there are $S$-equivariant maps $E'\rightarrow
E\rightarrow E''$, and hence augmentation preserving $RS$-maps 
$W'_*\rightarrow W_*\rightarrow W''_*$.  Thus it suffices to verify
the lower bound in the case when $W_*=W'_*$, and the upper bound in
the case when $W_*=W''_*$.  The lower bound follows from the Ku\"nneth
theorem, because 
$$\tot\Hom_S(W'_*,D(n,0)^*)\cong D(n,0)^* \cong C(n,0)^{*\otimes
\Omega}.$$

For the upper bound in the case when $\alpha$ generates a summand
$R/(n)$ for $n\neq 0$, we consider the type II spectral sequence,
$E_*^{*,*}$, for
$$E_0^{i,j}= \Hom_S(W''_i,D(n,0)^j).$$
Since we assumed that $W''_0$ is free of rank one, $E_0^{0,0}$ is
isomorphic to $R$, and some generator for $E_0^{0,0}$ is a cocycle
representing $\wre {\alpha'}$.  Since $E_0^{i,j}$ is zero for $j>0$, 
$E_\infty^{0,0}$ is a subgroup of $H^0(\tot E_0^{*,*})$.  Hence the
order of $\wre {\alpha'}$ is equal to the order of $E_\infty^{0,0}$.
Now as $RS$-module, $D(n,0)^{-1}$ is isomorphic to the permutation
module with basis $\Omega$.  (See the analysis at the end of
Section~1.)  The coboundary 
$$R\Omega = D(n,0)^{-1}\mapright{d} D(n,0)^0 = R$$ 
satisfies $d(\omega)= n$ for each $\omega \in \Omega\subseteq
R\Omega$.  Now since $W''_*$ is acyclic, the $E_1$-page of the
spectral sequence is isomorphic to the cohomology of the group $S$
with coefficients in $D(n,0)^*$.  More precisely, 
$$E_1^{i,j} = H^i(S,D(n,0)^j).$$
In particular, $E_1^{0,j}= (D(n,0)^j)^S$, the $S$-fixed points in
$D(n,0)^j$.  Thus $E_1^{0,0} = R$, and $E_1^{0,-1}$ is the free
$R$-module with basis the $S$-orbits in $\Omega$.  
It follows that 
$$E_2^{0,0} = E_1^{0,0}/\im(d':E_1^{0,-1}\rightarrow E_1^{0,0}) 
\cong R/(nl'),$$
where $l'$ is the h.c.f.\ of the lengths of the $S$-orbits in
$\Omega$.  This gives the required upper bound, since $E_\infty^{0,0}$
is a quotient of $E_2^{0,0}$. \qed

\remark Ther is an easy argument using the transfer which gives the
weaker upper bound $O(\wre \alpha)\leq |S|.O(\alpha)$.  

The most interesting case of the construction of $\wre\alpha$ is of
course the case in which $W_*$ is acyclic, and we shall concentrate on
that case from now on.  Proposition~3.2 gives an easy case in which
the lower bound given above is attained for $W_*$ acyclic, and hence
for all $W_*$.  

\proclaim Proposition 3.2.  With notation as in Theorem~3.1, if
$\alpha$ generates a summand of order $O(\alpha)= n$ with $n$ coprime
to $l'$, then $O(\wre\alpha)= O(\alpha)$.  

\proof By assumption, each of $\alpha$ and $l'\alpha$ generates a
summand of $H^*(C^*)$ of order $n$, so $O(\wre\alpha) = O\wre
{(l'\alpha)}$ and $O(\wre\alpha)$ divides $nl'$ by 3.1.  But  
$\wre {(l'\alpha)} = l'^l(\wre\alpha) $, and hence $O(\wre\alpha)$ is
coprime to $l'$.  \qed

In view of Proposition 3.2, it is reasonable to consider the problem
of the order of $\wre\alpha$ one prime at a time.  An easy transfer
argument shows that in this case it suffices to consider the case when
$S$ is a $p$-group.  In this case, and for $W_*$ acyclic, I know of no
case where the upper bound of Theorem~3.1 (or rather its $p$-part) is not
attained.  The only cases in which I have been able to prove this are
stated below.   

\proclaim Theorem 3.3. Let $S$ be a $p$-group of order $p^n$, and let
$S$ permute $\Omega$ freely, transitively.  Let $W_*$ be acyclic and
let  $\alpha$ generate a summand of $H^*(C^*)$ of order $p^r$.  Then 
$O(\wre\alpha) = p^{r+n}$ if either $S$ is cyclic, or $S$ is
elementary abelian.  

\proof By Theorem~3.1, it suffices to show this for some single choice
of $C^*$ and element $\alpha$ generating a summand of $H^*(C^*)$ of
even degree.  We choose $C^*$ to be a cochain complex for computing
the integral cohomology of some $p$-group $H$, and then use techniques
from cohomology of finite groups, including the Evens norm map [\evb] 
and a theorem of J.~F. Carlson [\carl] which we state below as Theorem~3.4.  

Recall (from for example [\evb]) that if $G$ is a finite group, with
subgroup $H$, then a choice of transversal $T$ to $H$ in $G$ gives
rise to an injective homomorphism 
$$\phi_T: G \mapright{} H\wreath \Sigma(G/H), $$ 
where $\Sigma(G/H)$ is the permutation group on the set $G/H$ of
cosets of $H$ in $G$.  
Furthermore, if $T'$ is another transversal, then $\phi_T$ and
$\phi_{T'}$ differ only by an inner automorphism of $H\wreath
\Sigma(G/H)$.  If $H$ is normal in $G$ then the map $G\rightarrow 
\Sigma(G/H)$ factors through $G/H$, and $G/H$ acts freely,
transitively on the cosets of~$H$.  If $\alpha$ is an element of $H^*(H)$ of
even degree, then the Evens norm $\norm^G_H(\alpha)\in H^*(G)$ is
defined to be $\phi_T^*(\wre\alpha)$ (for $W_*$ acyclic), and is
independent of $T$.  There is a double coset formula for the image of
$\norm^G_H(\alpha)$ in $H^*(H)$, which in the case when $H$ is normal
in $G$ is:  
$$\res^G_H\norm^G_H(\alpha) = \prod_{t\in T}c_t^*(\alpha),\eqno{(1)}$$
where $c_t^*$ is the map of $H^*(H)$ induced by conjugation by $t$.  
Our strategy for proving the theorem is to find a group $G$
expressible as an extension with quotient $S$ and kernel some suitable
$H$, and some $\alpha\in H^*(H)$ generating a summand of order $p^r$
such that $\norm^G_H(\alpha)$ has order $p^{n+r}$.  

In the case when $S= C_{p^n}$ is cyclic, we may take $G$ to be cyclic
of order $p^{n+r}$, so that $H$ is cyclic of order $p^r$.  Then 
$$H^*(H)= \Bbb Z[\alpha]/(p^r\alpha), \qquad 
H^*(G)= \Bbb Z[\alpha']/(p^{n+r}\alpha'), $$
where $\alpha$ and $\alpha'$ have degree two, 
and $\gamma\in H^{2i}(G)$ generates $H^{2i}(G)$ if and only if
$\res^G_H(\gamma)$ generates $H^{2i}(H)$.  
But by (1), $\res^G_H\norm^G_H(\alpha) = \alpha^{p^n}$, so
$\norm^G_H(\alpha)$ has order $p^{n+r}$ as required.  The case when
$S$ is non-cyclic of order four may be proved similarly, taking $G$
to the the generalized quaternion group of order $2^{r+2}$ expressed
as a central extension with quotient $S$ and cyclic kernel.  

The case $S= (C_p)^n$ for $n\neq 1$ and $(p,n)\neq (2,2)$ is more
complicated, because here it seems to be impossible to choose $G$ 
with $H$ a central subgroup.  For $p=2$, let $P$ be the dihedral group
of order eight, and for odd $p$, let $P$ be the (unique) non-abelian
group of order $p^3$ and exponent $p$.  In each case the centre of $P$
is cyclic of order $p$.  Now let $G$ be the central product of $n$
copies of $P$ and a cyclic group of order $p^r$.  The group $G$ has
the following presentation:  
$$G=\langle A_1,\ldots,A_n,B_1,\ldots,B_n,C\mid 
C^{p^r},\, A_i^p,\, B_i^p,\, [A_i,A_j],\quad$$
$$\qquad 
[B_i,B_j],\, [A_i,C],\, [B_i,C],\, [A_i,B_j]C^{-\delta(i,j)p^{r-1}}
\rangle \eqno{(2)} $$
where $\delta(i,j)$ is the Dirac $\delta$-function.  Let $Z$ be the
subgroup of $G$ generated by $C$.  Then $Z$ is the centre of $G$, 
$Z$ is cyclic of order $p^r$, and $G/Z$ is elementary abelian of
rank~$2n$.  Let $H$ be the subgroup of $G$ generated by $C$, and the
$A_i$'s.  Then $H$ is normal in $G$, and $H= Z\times H'$, where $H'$
is the elementary abelian group of rank $n$ generated by the $A_i$'s.
The quotient $G/H$ is elementary abelian of rank $n$, so is isomorphic
to $S$ as required.  

Since $Z$ is a direct summand of $H$, the map
$\res^H_Z:H^*(H)\rightarrow H^*(Z)$ is surjective.  Let $H^*(Z)= 
\Bbb Z[\gamma]/(p^r\gamma)$, and let $\alpha\in H^2(H)$ be such that 
$\res^H_Z(\alpha)=\gamma$.  The exponent of $H^*(H)$ is $p^r$ (by the
K\"unneth theorem), so any such $\alpha$ generates a summand of
$H^*(H)$.  Since $Z$ is central in $G$, $c_g^*$ acts trivially on
$H^*(Z)$ for any $g\in G$.  Applying the double coset formula given in
(1), it follows that 
$$\res^G_Z\norm^G_H(\alpha) = \prod_{t\in G/H} c_t^*(\gamma) =
\gamma^{p^n}. $$
The claim will follow if we can prove that any element of $H^*(G)$
whose image generates $H^{2i}(Z)$ for some $i>0$ has order $p^{n+r}$.
This is done in Lemma~3.5.  \qed

\proclaim Theorem 3.4. ([\carl]) Let $G$ be a finite group, and let
$x_1,\ldots, x_m$ be elements of $H^*(G)$ such that $H^*(G)$ is a
finite module for the subalgebra they generate.  Then the order $|G|$
of $G$ divides the product $O(x_1)\cdots O(x_m)$.  \qed

\proclaim Lemma 3.5. Let $G$ be the group with presentation (2) as
above, and $Z$ the subgroup generated by $C$.  Then any element of
$H^*(G)$ whose image under $\res^G_Z$ generates $H^{2i}(Z)$, for
some $i>0$, has order $p^{n+r}$.  

\proof Let $\alpha_1$ be an element as in the statement.  We shall
exhibit $\alpha_2,\ldots,\alpha_{n+1}\in H^*(G)$ of order $p$ such
that $H^*(G)$ is finite over the subalgebra generated by
$\alpha_1,\ldots,\alpha_{n+1}$.  Then by Theorem~3.4, $\alpha_1$ must
have order at least $p^{n+r}$.   The subgroup $H$ of $G$ has index 
$p^n$, and its (positive degree) cohomology has exponent $p^r$, so 
a transfer argument gives the other inequality.  Recall that work of
Quillen implies that for any $p$-group $K$, $H^*(K)$ is finite over a
subring $H'^*$ if and only if for every maximal elementary abelian
subgroup $E$ of $K$, $H^*(E)$ is finite over its subring
$\res^K_E(H'^*)$ [\quib].  

The maximal elementary abelian subgroups of $G$ have $p$-rank $n+1$,
and all contain the central subgroup of $Z$ order $p$.  (One way to
see this is to note that the subgroup of $G$ generated by the $A_i$'s
and $B_i$'s contains all elements of $G$ of order $p$, and this group
is an extraspecial group of order $p^{2n+1}$, whose elementary abelian
subgroups are discussed in for example [\bencarl].)  The quotient
$G/Z$ is elementary abelian of rank $2n$, and if $E$ is a maximal
elementary abelian subgroup of $G$, the image of $E$ in $G/Z$ is
elementary abelian of rank $n$.  

Recall that if $E$ is an elementary abelian group of rank $m$, then 
$H^2(E)\cong \Hom(E,\Bbb Q/\Bbb Z)$ is elementary abelian of the same
rank.  Let $H^{\prime *}$ be the subalgebra of $H^*(E)$ generated by $H^2$.
Elements of $H^*(E)$ of positive degree have exponent $p$, and 
$H^{\prime *}\otimes \fp$ is naturally isomorphic to the ring $\fp[E]$
of polynomial functions on $E$ viewed as an $\fp$-vector space.  Recall from 
for example [\beninv, Chap.~8], that the ring of invariants in $\fp[E]$
under the action of the full automorphism group of $E$ is a polynomial
algebra with generators $c_{m,0},\ldots,c_{m,m-1}$, where $c_{m,i}$
has degree $2(p^m-p^i)$.  (The $c_{m,j}$'s are known as the Dickson
invariants.)  Recall also that if $E'$ is a subgroup of
$E$ of rank $m'$, then the image of $c_{m,j}$ in $\fp[E']$ is zero 
for $j<m-m'$ and is a power of $c_{m',j-m+m'}$ otherwise.  

Now let $\gamma_1,\ldots,\gamma_n$ be the elements of $H^*(G/Z)$
corresponding to the Dickson invariants $c_{2n,n},\ldots,c_{2n,2n-1}$
in the subalgebra generated by $H^2(G/Z)$.  Let $E$ be a maximal
elementary abelian subgroup of $G$, and let $Z'$ be the order $p$
subgroup of $Z$.  Then $E/Z'$ is a subgroup of $G/Z$ of rank $n$, and 
the properties of the Dickson invariants stated above imply that 
$H^*(E/Z')$ is finite over the subalgebra generated by
$\res(\gamma_1),\ldots, \res(\gamma_n)$.  Let $\alpha_{i+1}$ be the
image of $\gamma_i$ in $H^*(G)$.  Now
$\res^G_E(\alpha_2),\ldots,\res^G_E(\alpha_{n+1})$ 
freely generate a polynomial subalgebra of $H^*(E)\otimes \fp$ and 
$\res^G_{Z'}(\alpha_i)= 0$ for $i>1$.  If $\alpha_1$ is an element of 
$H^*(G)$ with $\res^G_{Z'}(\alpha_1)\neq 0$, as in the statement, 
then it follows that $\res^G_E(\alpha_1),\ldots,\res^G_E(\alpha_{n+1})$
also freely generate a polynomial subalgebra of $H^*(E)\otimes \fp$,
and so $H^*(E)$ is finite over the subalgebra they generate.  
Thus $H^*(G)$ is finite over the algebra generated by
$\alpha_1,\ldots, \alpha_{n+1}$ by Quillen's theorem, and then 
$\alpha_1$ must have order at least $p^{n+r}$ by Carlson's
Theorem~3.4.  \qed

\remark The proof of Lemma~3.5 is based on Carlson's method for
computing the exponent of the cohomology of the extraspecial groups
[\carl].  Carlson does not use Dickson invariants to construct
elements analogous to the $\alpha_i$'s for $i>1$, but gives an 
existence proof for such elements via algebraic geometry.

\beginsection 4. INTEGRAL COHOMOLOGY FOR THE CYCLIC GROUP OF ORDER $p$

In this section we concentrate on the study of the integral cohomology
of $X^p\times_{C_p}E$, where $C_p$ acts by freely permuting the
factors of $X^p$ and $E$ is contractible.  As above, the results apply
more generally to the cohomology of $\Hom_{C_p}(W_*,C^{*\otimes p})$
for any finite type cochain complex $C^*$ of free abelian groups.  
Evens and Kahn
obtained partial results for this case [\evk] and we give few
details for that part of our calculation which is a repeat of
theirs.  The main new idea here is the use of Lemma 1.1, which
enables us to complete the calculation of 
$H^*(X^p\times_{C_p}E)$ and to give simpler proofs of some of 
the results in Evens and Kahn's paper [\evk].  

As in Section 1, let $C^*$ be a cochain complex of finitely generated
free abelian groups, for example, the cellular cochain complex of a
CW-complex of finite type.  Let $E$
be a contractible $C_p$-free $C_p$-CW-complex, and let $W_*$ be
the cellular chain complex for $E$.  Recall from Lemma~1.4 that 
for finding the cohomology of $\Hom_{C_p}(W_*,C^{*\otimes p})$, 
the complex $C^*$ may be replaced by any homotopy equivalent complex.  
As in Section~1, we replace $C^*$ by a complex $C'^*$ consisting of a direct
sum of pieces of the form $C(n,i)^*$, in bijective correspondence with
the summands of $H^i(C^*)$ (in some fixed splitting) isomorphic to $\zz/(n)$.  
(Recall that $C(n,i)^*$ is a complex which has at most two non-zero
groups, each of rank one, and that $H^*(C(n,i)^*)$ is isomorphic to
$\zz/(n)$ concentrated in degree $i$.)

As in Section 1, the complex of abelian groups, $C'^{*\otimes p}$ 
splits as a
direct sum of \lq cubes' of the form $C(n_1,i_1)^*\otimes\cdots
\otimes C(n_p,i_p)^*$.  The action of $C_p$ permutes the cubes freely,
except for those cubes whose sides all correspond to the same
cyclic summand of $H^*(C^*)$.  The action of $C_p$ on these cubes 
is by a cyclic permutation of the $p$ distinct axes.  Considering the
cube as embedded in $\Bbb R^p$, the matrix for the action is a
permutation matrix of order $p$.  More care must be taken for $p=2$
than for $p$ odd, because a $2\times 2$ permutation matrix of order
two has determinant $-1$.

The cubes permuted freely by $C_p$ cause no problem.  Indeed, if 
$D^*$ is any cochain complex of abelian groups, then the
Eckmann-Shapiro lemma [\bro,\ben] shows that
$$H^*\tot\Hom_{C_p}(W_*,D^*\otimes \Bbb ZC_p)\cong H^*(D^*).$$
Moreover, each of the two spectral sequences
arising from viewing $\Hom_{C_p}(W_*,D^*\otimes \Bbb ZC_p)$ as a
double cochain complex has $E_2^{i,0}\cong H^i(D^*)$ and
$E_2^{i,j}=0$ for $j\neq 0$.  We can of course find the
cohomology of $C(n_1,i_1)\otimes\cdots\otimes C(n_p,i_p)$ using
the Kunneth theorem.  If $i=i_1+\cdots+i_p$ and exactly $r$ of 
$n_1,\ldots,n_p$ are nonzero, then $H^{i-j}(C(n_1,i_1)\otimes 
\cdots C(n_p,i_p))$ is isomorphic to a sum of ${r-1}\choose 
j$ copies of $\Bbb Z/(n_1,\ldots,n_p)$.  

Similarly, the cubes of the form $C(0,i)^{\otimes p}$ where
$C_p$ acts by permuting the factors are easy to handle.  The
cochain complex $C(0,i)^{\otimes p}$ consists of a single 
$\Bbb ZC_p$-module of $\Bbb Z$-rank one in degree $pi$.  This is
the trivial $\Bbb ZC_p$-module except when $p=2$ and $i$ is odd,
in which case it is $\hat {\Bbb Z}$, the module on which a
generator for $C_2$ acts as multiplication by $-1$.  Thus each
such cube contributes a summand to the spectral sequence of the 
form $H^*(C_p;\Bbb Z)$ (resp. $H^*(C_2;\hat{\Bbb Z})$ if $p=2$ and 
$i$ is odd) concentrated in $E_2^{*,pi}$.  

As in Section~3, let $D(n,i)^*=C(n,i)^{\otimes p}$ for 
$n\geq 0$ be a complex
of $\Bbb ZC_p$-modules where $C_p$ acts by permuting the factors
(with a sign if $p=2$ and $i$ is odd).  The only contributions 
to $H^*(X^p\times_{C_p}E)$ not accounted for by the above
remarks come from summands of the double cochain complex of the 
form $\Hom_{C_p}(W_*,D(n,i))$, for $n>1$.  It is easy to see that the
module $D(n,i)^j$ is the zero module unless $p(i-1)\leq j\leq
pi$, and that $D(n,i)^j$ is $\Bbb ZC_p$-free of rank
$1/p{p\choose j}$ if $p(i-1)<j<pi$.  In the case when $p$ is
odd, $D(n,i)^{pi}$ and $D(n,i)^{p(i-1)}$ are the trivial 
$\Bbb ZC_p$-module $\Bbb Z$.  When $p=2$, $D(n,i)^{2i}$ is 
isomorphic to $\Bbb Z$ (resp.\ $\hat{\Bbb Z}$) and
$D(n,i)^{2(i-1)}$ is isomorphic to $\hat {\Bbb Z}$ (resp.\ $\Bbb
Z$) if $i$ is even (resp.\ odd).  To describe the differential
in $D(n,i)^*$ completely we would have to choose an explicit 
generating set for each $D(n,i)^j$, and we shall not do this.
All that we shall require in the sequel is that $D(1,i)^*$ is 
exact (because $C(1,i)^*$ is exact), and that if we view
$C(1,i)^*$ and $C(n,i)^*$ as consisting of the same groups but with
different maps, then the differential in $D(n,i)^*$ is $n$-times
that in $D(1,i)^*$.  Note that for $p$ odd, each $D(n,i)$ is
isomorphic to $D(n,0)$ shifted in degree by $pi$, and for $p=2$, 
$D(n,2i)$ is isomorphic to $D(n,0)$ shifted by $4i$ and 
$D(n,2i+1)$ is isomorphic to $D(n,1)$ shifted by $4i$.

\proclaim Theorem 4.1.  For $p$ an odd prime, let $E_0^{i,j}
=\Hom_{C_p}(W_i,D(n,0)^j)$, and let
$E_*^{*,*}$ be the corresponding type II spectral sequence.  Then for
$i>0$, $E_2^{i,j}$ is zero except that $E_2^{2i,0}$ and $E_2^{2i,-p}$
are cyclic of order $p$.  Define a function $g(j)$ as follows: 
$$g(j)=\cases{{1/p}((p-1)(-1)^j+{p-1\choose j})&for $0\leq j\leq p-1$,\cr
0&otherwise.\cr}$$
Then $E_2^{0,-j}$ is an extension of $\Bbb Z/p$ by $(\Bbb Z/n)^{g(j)}$ 
for $j=0,2,4,\ldots,p-3$, nonsplit if $p$ divides $n$, and 
for other $j$, $E_2^{0,-j}$ is isomorphic to $(\Bbb Z/n)^{g(j)}$.  
If $p$ divides $n$ then the spectral sequence collapses at $E_2$.  
If $p$ does not divide $n$ then the non-zero higher differentials 
are $d_3,d_5,\ldots,d_p$, and in this case $E_\infty^{i,j}=0$ 
for $i>0$ and $E_\infty^{0,-j}\cong (\Bbb Z/n)^{g(j)}$.  

\remark The reader may find it helpful to consult figures 1~and~2, 
which illustrate cases of the above statement.  

\proof For this spectral sequence $E_1^{i,j}$ is isomorphic to 
$H^i(C_p;D(n,0)^j)$.  Since $D(n,0)^{-j}$ is the trivial module 
$\Bbb Z$ for $j=0$ or $j=p$, free of rank ${1/p}{p\choose j}$ 
for $0<j<p$ and zero for other $j$, we see that $E_1^{i,j}$ is
trivial for $i>0$ except that $E_1^{2i,0}\cong E_1^{2i,-p}
\cong \Bbb Z/p$.  Also $E_1^{0,-j}$ is free abelian of rank 
$f(j)$, where $f$ is defined as follows:  
$$f(j)=\cases{1&for $j=0$ or $j=p$,\cr
1/p{p\choose j}&for $0<j<p$,\cr
0&otherwise.\cr}$$
Note that the only non-zero groups on the $E_1$-page occur on the
three line segments $j=0$ and $i\geq 0$, $j=-p$ and $i\geq 0$, 
$i=0$ and $-p\leq j\leq 0$.  The shape of the $E_1$-page implies 
that $E_2^{i,j}=E_1^{i,j}$ for $i>0$, and that the only possibly 
non-zero differentials after $d_1$ are the following:  
$$\eqalign{d_3:E_3^{2,-p}&\rightarrow E_3^{0,3-p},\cr
d_5:E_5^{4,-p}&\rightarrow E_5^{0,5-p},\cr
&\vdots\cr
d_{p-2}:E_{p-2}^{p-3,-p}&\rightarrow E_{p-2}^{0,-2},\cr
\hbox{ and }\quad d_p: E_p^{p-1+2i,-p}&\rightarrow E_p^{0,2i}
\quad\hbox{ for all $i\geq 0$.}\cr}$$
To determine the groups $E_2^{0,j}$ and the higher differentials we
shall first consider the case $n=1$ and then apply Lemma~1.1.  

The complex $D(1,0)^*$ is exact, which implies that for $n=1$ the 
total complex of $E_0^{*,*}$ is exact, and hence $E_\infty^{i,j}=0$
for all $i$ and $j$.  We also know the isomorphism type of
$E_2^{i,j}$ for all $i$ and $j$ except for the cases when $i=0 $ 
and $-p\leq j\leq 0$.  Since each of the possibly non-trivial groups
on the $E_2$-page is involved in at most one possibly non-trivial 
higher differential, it follows that in the case $n=1$,	 all the 
possibly non-zero differentials listed above must be isomorphisms,
and that all groups $E_2^{0,j}$ except those such that $E_r^{0,j}$ 
appears in the above list must be trivial.  This completes the 
proof in the case $n=1$.  

\midinsert
$$\spreaddiagramrows{-1.5pc}
\spreaddiagramcolumns{-1.5pc} 
\diagram
\bu&\mi&\bu&\mi&\bu&\mi&\bu&\mi&\bu\\
\mi&\mi&\mi&\mi&\mi&\mi&\mi&\mi&\mi&\mi\\
\bu&\mi&\mi&\mi&\mi&\mi&\mi&\mi&\mi&\mi\\
\mi&\mi&\mi&\mi&\mi&\mi&\mi&\mi&\mi&\mi&\mi\\
\bu&\mi&\mi&\mi&\mi&\mi&\mi&\mi&\mi&\mi&\mi\\
\mi&\mi&\mi&\mi&\mi&\mi&\mi&\mi&\mi&\mi&\mi&\mi\\
\mi&\mi&\mi&\mi&\mi&\mi&\mi&\mi&\mi&\mi&\mi&\mi\\
\mi&\mi&\bu\ar[uuull]_{d_3}&\mi&\bu\ar[uuuuullll]_{d_5} 
&\mi&\bu\ar[uuuuuuullllll]_{d_7}
&\mi&\bu\ar[uuuuuuullllll]_{d_7}
&\mi&\bu\ar[uuuuuuullllll]_{d_7}
&\mi&\bu\ar[uuuuuuullllll]_{d_7}
\enddiagram 
$$

\centerline{Fig.\ 1. The $E_2$-page and higher differentials for 
the spectral sequence of 4.1,} 
\centerline{when $p=7$ and $n=1$.  \lq$\bullet$' denotes a non-zero entry.}  
\endinsert

\midinsert
$$
\vbox{\tabskip=0pt \offinterlineskip
\def\tablerule{\noalign{\hrule}}
\halign{
\bigstrut#&\strut#\tabskip=1em plus2em& 
\hfil$#$\hfil&\strut#&
\hfil$#$\hfil&\strut#&
\hfil$#$\hfil&\strut#&
\hfil$#$\hfil&\strut#&
\hfil$#$\hfil&\strut#&
\hfil$#$\hfil&\strut#&
\hfil$#$\hfil&\strut#&
\hfil$#$\hfil\tabskip=0pt\crcr
\tablerule
&\vrule&\ze/7n&\vrule&\mj&\vrule&\ze/7&\vrule&\mj&\vrule&\ze/7&\vrule
&\mj&\vrule&\ze/7&\vrule&\cdots\ \ \cr
\tablerule
&\vrule&\mj&\vrule&\mj&\vrule&\mj&\vrule&\mj&\vrule&\mj&\vrule
&\mj&\vrule&\mj&\vrule&\cdots\ \ \cr
\tablerule
&\vrule&\ze/7n\oplus(\ze/n)^2
&\vrule&\mj&\vrule&\mj&\vrule&\mj&\vrule&\mj&\vrule
&\mj&\vrule&\mj&\vrule&\cdots\ \ \cr
\tablerule
&\vrule&(\ze/n)^2
&\vrule&\mj&\vrule&\mj&\vrule&\mj&\vrule&\mj&\vrule
&\mj&\vrule&\mj&\vrule&\cdots\ \ \cr
\tablerule
&\vrule&\ze/7n\oplus(\ze/n)^2
&\vrule&\mj&\vrule&\mj&\vrule&\mj&\vrule&\mj&\vrule
&\mj&\vrule&\mj&\vrule&\cdots\ \ \cr
\tablerule
&\vrule&\mj&\vrule&\mj&\vrule&\mj&\vrule&\mj&\vrule&\mj&\vrule
&\mj&\vrule&\mj&\vrule&\cdots\ \ \cr
\tablerule
&\vrule&\ze/n
&\vrule&\mj&\vrule&\mj&\vrule&\mj&\vrule&\mj&\vrule
&\mj&\vrule&\mj&\vrule&\cdots\ \ \cr
\tablerule
&\vrule&\mj&\vrule&\mj&\vrule&\ze/7&\vrule&\mj&\vrule&\ze/7&\vrule
&\mj&\vrule&\ze/7&\vrule&\cdots\ \ \cr
\tablerule}}
$$

\centerline{Fig.\ 2. The $E_2$-page of the spectral 
sequence of 4.1 for $p=7$.}
\endinsert

Let $Z^j$ and $B^j$ stand for the cycles and boundaries respectively
in $E_1^{0,-j}$ for the case $n=1$.  Then $Z^j$ and $B^j$ are free 
abelian of the same rank, and $Z^j=B^j$  except that $Z^j/B^j$ has
order $p$ for $j=0,2,4,\ldots,p-3$.  Let $g'(j)$ stand for the 
rank of $Z^j$ or $B^j$.  By Lemma~1.1., the group of
cycles for $d_1$ in $E_1^{0,-j}$ for general $n$ is equal to $Z^j$,
while the group of boundaries is equal to $nB^j$.  It follows that
for general $n$, $E_2^{0,-j}$ is the natural extension of $Z^j/B^j$ by 
$B^j/nB^j \cong (\Bbb Z/n)^{g'(j)}$, which is nonsplit whenever
possible.  Thus to verify
the claimed description of the $E_2$-page it suffices to show that 
$g'(j)=g(j)$.  For this, note that the short exact sequence 
$$0\rightarrow Z^j\rightarrow E_1^{0,-j}\rightarrow 
B^{j-1}\rightarrow 0$$ 
implies that 
$g'(j)+g'(j-1)=f(j)$.  Define polynomials $F(t)$, $G(t)$ in a formal
variable $t$ by 
$$F(t)=\sum_j f(j)t^j,\qquad G(t)=\sum_jg'(j)t^j,$$
and note that 
$$F(t)={p-1\over p}(1+t^p) + {1\over p}(1+t)^p.$$
The relation given between $f$ and $g'$ implies that
$(1+t)G(t)=F(t)$, from which it is easy to verify that $g'=g$.  

It now remains only to check the given description of the higher
differentials in the spectral sequence for general $n$.  Once more we
invoke Lemma~1.1.  For $i>0$, let $x\in E_0^{2i,-p}$ be an element 
representing a generator for $E_1^{2i,-p}$, and let $y\in
E_0^{2i-1-p}$ (resp.\ $y\in E_0^{2i-p+1,0}$ if $2i>p$) be an element
representing the image of $x$ under $d_{2i+1}$ (resp.\ $d_p$) in the
spectral sequence for $n=1$.  If $p$ divides $n$, then $n^jy$ will
represent zero in $E_2$ for $j\geq 2$, and so in this case the higher
differentials are trivial.  If on the other hand $p$ does not divide
$n$ then $n^jy$ will represent an element of order $p$ in $E_2$ for 
$j\geq 1$, and so again the higher differentials are as claimed.  
\qed

\def\pra{\par}
\proclaim Theorem 4.1$'$.  For $p=2$, let $E_0^{*,*}$ be the double 
cochain complex $\Hom_{C_2}(W_*,D(n,0)^*)$.  Then in the
corresponding type II spectral sequence, $E_2^{i,j}=0$ except that 
$E_2^{2i+2,0}\cong E_2^{2i+1,-2}\cong \Bbb Z/2$ for all $i\geq 0$,
and $E_2^{0,0}\cong \Bbb Z/(2n)$.  If $n$ is even the spectral
sequence collapses at $E_2$.  If $n$ is odd, the spectral sequence 
collapses at $E_3$, and $E_3^{i,j}=0$ except that $E_3^{0,0}\cong 
\Bbb Z/n$.  
\pra 
Let ${E'}_0^{*,*}$ be the double cochain complex
$\Hom_{C_2}(W_*,D(n,1)^*)$.  Then in the corresponding type II
spectral sequence $E_2^{i,j}=0$ except that ${E'}_2^{2i+1,2}\cong 
{E'}_2^{2i+2,0}\cong \Bbb Z/2$ for all $i\geq 0$, and
${E'}_2^{0,1}\cong \Bbb Z/n$.  If $n$ is even the spectral sequence
collapses at $E_2$.  If $n$ is odd the spectral sequence collapses 
at $E_3$ and ${E'}_3^{i,j}=0$ except that ${E'}_3^{0,1}\cong \Bbb
Z/n$.  

\proof Similar to the proof of Theorem~4.1, and easier in spite of the
extra complication introduced by the second action of $C_2$ on 
$\Bbb Z$.  
\qed

\proclaim Theorem 4.2.  Let $p$ be an odd prime and let $W_*$ be an 
acyclic complex of free $\zz C_p$-modules.  Let $C^*$ be a 
cochain complex of finitely generated free abelian groups, 
and choose some splitting 
$$H^*(C^*)\cong \bigoplus_{ a\in A}H( a),$$
where $H( a)$ is a summand of $H^{d( a)}(C^*)$, and is
isomorphic to $\Bbb Z/(n( a))$, where (for simplicity and 
without loss of generality) each $n( a)$ is either 0, a
power of $p$, or a positive integer coprime to $p$.  Let $C_p$ act on
the set $A^p$ by permuting the factors, so that elements of the 
form $( a_1,\ldots, a_p)$ are in free orbits provided that 
not all the $ a_i$ are equal, and elements of the form 
$( a,\ldots, a)$ are fixed.  Then $H^*\tot\Hom_{C_p}(W_*,C^{*\otimes p})$
is a direct sum of the following summands.  
\pra 
\item{a)} For each free $C_p$-orbit in $A^p$ with orbit representative 
$( a_1,\ldots, a_p)$, and for each $j$, a direct sum of 
${r\choose j-1}$ copies of $\Bbb Z/(n( a_1),\ldots,n( a_p))$ 
in degree $d( a_1)+\cdots+d( a_p)-j$.  Here $r$ is the number
of $ a_i$'s such that $n( a_i)$ is non-zero.  
\pra
\item{b)} For each $ a$ such that $n( a)=0$, one copy of $\Bbb Z$ in 
degree $pd( a)$, and for  each $i\geq 0$, one copy of 
$\Bbb Z/(p)$ in degree $pd( a)+2+2i$.  
\pra
\item{c)} For each $ a$ such that $n( a)$ is non-zero and coprime to
$p$ and for each $j$, a direct sum of $g(j)$ copies of $\Bbb
Z/(n a)$ in degree $pd( a)-j$, where $g(j)$ is as defined in 
the statement of Lemma 2.2.
\pra
\item{d)} For each $ a$ such that $n( a)$ is a power of $p$, for 
$j=0,2,\ldots,p-3$, a direct sum of $g(j)-1$ copies of $\Bbb
Z/(n( a))$ and one copy of $\Bbb Z/(pn( a))$ in degree
$pd( a)-j$; for other $j$ a sum of $g(j)$ copies of 
$\Bbb Z/(n( a))$ in degree $pd( a)-j$; and for each 
$i\geq 0$ one copy of $\Bbb Z/(p)$ in degree $pd( a)+2+2i$ and 
one copy of $\Bbb Z/(p)$ in degree $p(d( a)-1)+2+2i$.  

\proof The fact that we can split the cohomology as a 
direct sum of contributions of the types described, and the analyses
of cases a) and b), were proved earlier.  
Cases c) and d) are descriptions of the cohomology of the total 
complex of $E_0^{*,*}=\Hom_{C_p}(W_*,D(n(\alpha),0))$, and so follow
from Theorem~4.1 except that in case d) we have to show that the
extension with kernel $E_\infty^{0,-p+2i}$ and quotient
$E_\infty^{2i,-p}$ representing $H^{2i-p}$ is split for $1\leq i\leq 
(p-1)/2$.  For $i=1$ this is obvious, because $g(p-2)=0$ and hence 
$E_\infty^{0,2-p}=0$.  For other $i$ we use (for the first time) the
graded $H^*(C_p;\Bbb Z)$-module structure of the spectral sequence,
which is a filtration of the graded $H^*(C_p)$-module structure on 
$H^*\Hom_{C_p}(W_*,D(n,0)^*)$.  It is easy to see that the product of
a generator for $E_\infty^{2,-p}$ and a generator for
$H^{2i}(C_p;\Bbb Z)$ is a generator for $E_\infty^{2+2i,-p}$ and
therefore that there is an element of $H^*(E_0^{*,*})$ of order $p$
yielding a generator for $E_\infty^{2+2i,-p}$, and so the extension
is split.  

As an alternative we may solve the extension problem by 
calculating $H^*(\tot E_0^{*,*}\otimes \Bbb F_p)$, which determines the 
number of cyclic summands of $H^*(\tot E_0^{*,*})$ by the universal
coefficient theorem.  When $p$ divides $n$, the differential in 
$D(n,0)\otimes \Bbb F_p$ is trivial, and so the type II spectral
sequence for $H^*(\tot E_0^{*,*}\otimes \Bbb F_p)$ collapses at the 
$E_1$-page, which makes this calculation easy.  We leave the 
details to the interested reader.  
\qed

\proclaim Theorem 4.2$'$.  Let $p=2$, and let $W_*$ be an acyclic 
complex of free $\zz C_2$-modules.  
Let $C^*$ be a cochain complex of finitely generated free abelian
groups and fix a splitting 
$$H^*(C^*)\cong\bigoplus_{ a\in A}H( a),$$
as in the statement of Theorem 4.2.  Then $H^*\Hom_{C_2}(W_*,C^{*\otimes
2})$ is a direct sum of the following summands.  
\pra
\item{a)} For each $ a\neq  a'\in A$, one summand $\Bbb
Z/(n( a),n( a'))$ in degree $d( a)+d( a')$, and if 
both $n( a)$ and $n( a')$ are nonzero, one summand $\Bbb
Z/(n( a),n( a'))$ in degree $d( a)+d( a')-1$.  
\pra
\item{b)} For each $ a$ with $n( a)=0$ and $d( a)$ even, one 
summand $\Bbb Z$ in degree $2d( a)$ and one summand $\Bbb Z/2$ in
each degree $2d( a)+2+2i$.  
\pra
\item{c)} For each $ a$ with $n( a)=0$ and $d( a)$ odd, one 
summand $\Bbb Z/2$ in each degree $2d( a)+1+2i$.  
\pra
\item{d)} For each $ a$ with $n( a)$ odd and $d( a)$ even, 
one copy of $\Bbb Z/(n( a))$ in degree $2d( a)$.  
\pra
\item{e)} For each $ a$ with $n( a)$ odd and $d( a)$ odd, 
one copy of $\Bbb Z/(n( a))$ in degree $2d( a)-1$.  
\pra
\item{f)} For each $ a$ with $n( a)$ a (strictly positive) power 
of 2 and $d( a)$ even, one copy of $\Bbb Z/(2n( a))$ in 
degree $2d( a)$ and one copy of $\Bbb Z/2$ in degree
$2d( a)-1$, and in each degree of the form $2d( a)+1+i$.  
\pra
\item{g)} For each $ a$ with $n( a)$ a (strictly positive) power of 
2 and $d( a)$ odd, one copy of $\Bbb Z/(n( a))$ in degree 
$2d( a)-1$ and once copy of $\Bbb Z/2$ in each degree
$2d( a)+i$.  

\proof This follows from the preamble together with Theorem~4.1$'$.  
Note that this is simpler than Theorem~4.2 in that there are no extension 
problems that need be resolved.  
\qed

\beginsection 5. THE SPECTRAL SEQUENCE FOR THE CYCLIC GROUP OF ORDER $p$

In the previous section we computed the integral cohomology of 
$X^p\times_{C_p}EC_p$ for any finite type CW-complex $X$, by
replacing the double cochain complex forming the $E_0$-page of the 
Cartan-Leray spectral sequence with a direct sum of simpler
complexes.  In this section we solve the associated type I spectral
sequences, and hence describe the differentials and extension
problems in the Cartan-Leray spectral sequence for
$X^p\times_{C_p}EC_p$.  

For $C^*$ a cochain complex of finitely generated free abelian groups, 
fix a splitting of $H^*(C^*)$ as a direct sum of cyclic groups, 
$$H^*(C^*)=\bigoplus_{a\in A}H(a)$$
as in the statement of Theorem~4.2.  
As in the preamble to Section~4, choose a cochain complex ${C'}^*$
splitting as a direct sum of subcomplexes $C'(a)^*$ (indexed by
the same set $A$), and a homotopy equivalence $f$ from ${C'}^*$ to 
$C^*$, the cochains on $X$, such that the image of $H^*(C'(a)^*)$ 
under $f^*$ is $H(a)\subseteq H^*(C^*)$.  As in Section~4, let 
$W_*$ be the chain complex for $EC_p$.  By Lemma~1.4, the double 
complexes 
$$E_0^{*,*}=\Hom_{C_p}(W_*,((C^*)^{\otimes p})^*)
\quad\hbox{and}\quad
{E'}_0^{*,*}=\Hom_{C_p}(W_*,(({C'}^*)^{\otimes p})^*)$$
give rise to isomorphic type I spectral sequences, from the
$E_1$-page onwards.  If $C^*$ is the cellular cochain complex of a
finite type CW-complex $X$, then 
the type I spectral sequence for $E_0^{*,*}$ is
the Cartan-Leray spectral sequence for $X^p\times_{C_p}EC_p$.  

We have already seen that ${E'}_0^{*,*}$ splits as a direct sum of 
subcomplexes indexed by the $C_p$-orbits in $A^p$.  In the preamble 
to Section~4 we showed that the type I spectral sequences
corresponding to free $C_p$-orbits in $A^p$ have $E_2$-page
concentrated in the column ${E'}_2^{0,*}$, so collapse at $E_2$ and
give rise to no extension problems.  We also showed that the trivial 
$C_p$-orbits of the form $(a,\ldots,a)$ such that
$H(a)$ is infinite cyclic give rise to double complexes with a
single nonzero row, and hence the type I spectral sequences for such
orbits collapse at $E_2$ and give rise to no extension problems.  
Thus the Cartan-Leray spectral 
sequence for $X^p\times_{C_p}EC_p$ splits from $E_2$ onwards as 
a direct sum of the following:  

\item{a)} Various pieces concentrated in $E_*^{0,*}$.

\item{b)} For each $a\in A$ such that $H(a)\subseteq H^i(X)$ 
is infinite cyclic, a piece concentrated in the row $E_*^{*,pi}$.  

\item{c)} For each $\alpha$ such that $H(a)\subseteq H^i(X)$ is cyclic
of order $n$, a copy of the type I spectral sequence for the double
complex $\Hom_{C_p}(W_*,D(n,i)^*)$ defined in the introduction to 
Section~4.  

Thus to solve the Cartan-Leray spectral sequence, it
suffices to solve the type I spectral sequences for the double
complexes of Theorems 4.1 and 4.1$'$.  

\proclaim Theorem 5.1.  Let $p$ be an odd prime, and let 
$E_0^{*,*}$ be the double cochain complex
of Theorem~4.1, and let $E_*^{*,*}$ be the corresponding type I
spectral sequence.  If $p$ does not divide $n$, then $E_2^{i,j}$ 
is as follows, where the function $g$ is as defined in Theorem~4.1:  
$$E_2^{i,j}=\cases{(\Bbb Z/n)^{\oplus g(-j)}&for $i=0$,\cr
0&for $i>0$.}$$ 
In this case the spectral sequence clearly collapses and gives rise
to no extension problems.  
If $p$ divides $n$, the $E_2$-page is as follows:  
$$E_2^{i,j}=\cases{(\Bbb Z/n)^{\oplus g(-j)}\oplus \Bbb Z/p&for 
$i=0$, $j$ odd, $0>j>1-p$,\cr
(\Bbb Z/n)^{\oplus g(-j)}&for $i=0$, $j$ not as above,\cr
\relax\Bbb Z/p&for $i>0$, $0\geq j\geq 1-p$,\cr 
0&otherwise.}$$
In this case the spectral sequence collapses at $E_3$, and 
the $E_3$-page is as follows:  
$$E_3^{i,j}=\cases{(\Bbb Z/n)^{\oplus g(-j)}&if $i=0$,\cr
\relax\Bbb Z/p&if either $j=0$, $i>0$ and $i$ even,\cr 
&\qquad or $j=1-p$ and $i$ odd,\cr
&\qquad or $i=1$ and $j=-1,-3,\ldots,2-p$,\cr
0&otherwise.}$$
The only non trivial extensions in reassembling $H^*\tot E_0^{*,*}$ 
from $E_\infty^{*,*}$ are that the extension with kernel 
$E_3^{1,-j}$ and quotient $E_3^{0,1-j}$ is non-split for
$j=1,3,\ldots,p-2$.  

\remark Figure~3 illustrates the case of the above statement when
$p=7$ and $n$ is a multiple of $p$.  In the figure, circles indicate
entries isomorphic to $\Bbb Z/7$, and squares indicate other non-zero
entries.  Entries in black remain non-zero in $E_3=E_\infty$.  Arrows
represent non-zero $d_2$'s, and double lines represent non-split
extensions.  

\midinsert
$$\spreaddiagramrows{-1.3pc}
\spreaddiagramcolumns{-1.3pc} 
\diagram
\sq\ar@{=}[dr]&\ci\ar[drr]&\bu&\ci\ar[drr]&\bu&
\ci\ar[drr]&\bu&\ci\ar[drr]&\bu&\ci&\bu&\cdots\\
\ci\ar[drr]&\bu&\ci\ar[drr]&\ci&\ci\ar[drr]&\ci&
\ci\ar[drr]&\ci&\ci\ar[drr]&\ci&\ci&\cdots\\
\sq\ar@{=}[dr]&\ci\ar[drr]&\ci&\ci\ar[drr]&\ci&
\ci\ar[drr]&\ci&\ci\ar[drr]&\ci&\ci&\ci&\cdots\\
\sq\ar[drr]&\bu&\ci\ar[drr]&\ci&\ci\ar[drr]&\ci&
\ci\ar[drr]&\ci&\ci\ar[drr]&\ci&\ci&\cdots\\
\sq\ar@{=}[dr]&\ci\ar[drr]&\ci&\ci\ar[drr]&\ci
&\ci\ar[drr]&\ci&\ci\ar[drr]&\ci&\ci&\ci&\cdots\\
\ci\ar[drr]&\bu&\ci\ar[drr]&\ci&\ci\ar[drr]&\ci&
\ci\ar[drr]&\ci&\ci\ar[drr]&\ci&\ci&\cdots\\
\sq&\bu&\ci&\bu&\ci&\bu&\ci&\bu&\ci&\bu&\ci&\cdots
\enddiagram
$$
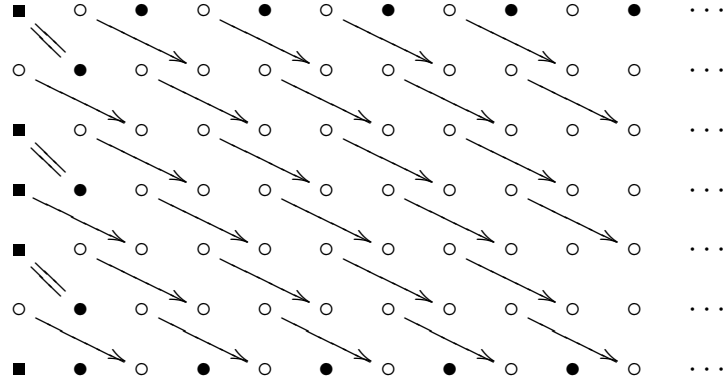
\centerline{Fig.\ 3. The $E_2$-page of the spectral sequence of 5.1,
for $p=7$ and $p | n$.}
%\centerline{Circles denote entries isomorphic to $\Bbb Z/7$, squares
%other non-zero entries.}  
%\centerline{Entries in black are non-zero in $E_3 = E_\infty$.  Arrows
%represent $d_2$} 
%\centerline{and double lines represent non-split extensions.}
\endinsert

\proof Recall that $E_0^{i,j}=\Hom_{C_p}(W_i,D(n,0)^j)$, and so 
$E_1^{i,j}=\Hom_{C_p}(W_i,H^j(D(n,0)^*))$ because $W_i$ is free, and
then $E_2^{i,j}=H^i(C_p;H^j(D(n,0)^*))$.  The complex $D(1,0)^*$ is
exact, and if we let $Z^j$ stand for the cycles in $D(1,0)^{-j}$, 
then for any $n$, 
$$H^{-j}(D(n,0)^*)=Z^j/nZ^j.$$
Since $Z^j$ is free abelian, the cochain complex
$\Hom_{C_p}(W_*,Z^j)$ is also free abelian, and we have the following
isomorphism.  
$$\Hom_{C_p}(W_*,Z^j/nZ^j)\cong \Hom_{C_p}(W_*,Z^j)/n\Hom_{C_p}(W_*,Z^j)$$
Hence we may compute $E_2^{i,-j}=H^i(C_p;Z^j/nZ^j)$ by first
computing $H^*(C_p;Z^j)$ and then applying a universal coefficient
theorem.  

Now $Z^0=D(1,0)^0=\Bbb Z$, $Z^j=0$ unless $0\leq j\leq p-1$, 
and for $0 <j\leq p-1$ there is a short exact sequence of 
$C_p$-modules
$$0\rightarrow Z^j\rightarrow F_j\rightarrow Z^{j-1}
\rightarrow 0$$
where $F_j$ is a free module of rank ${1\over p}{p\choose j}$.
Taking cohomology we obtain for $0<j<p$ the following exact
sequences.  
$$0\rightarrow H^0(C_p;Z^j)\rightarrow H^0(C_p;F_j)
\rightarrow H^0(C_p;Z^{j-1})\rightarrow H^1(C_p;Z^j)\rightarrow 0$$
$$0\rightarrow H^{i+1}(C_p;Z^{j-1})\rightarrow H^{i+2}(C_p;Z^j)
\rightarrow 0$$
Using the second of these and the periodicity of the cohomology of 
$C_p$ with arbitrary coefficients, it follows easily that for $i>0$ 
and $0\leq j\leq p-1$, 
$$H^i(C_p;Z^j)=\cases{\relax\Bbb Z/p&for $i+j$ even,\cr
0&for $i+j$ odd.}$$

For each $j$, $H^0(C_p,Z^j)$ is a free abelian group, so we need only
find its rank.  From the first of the two exact sequences we obtain 
for $0<j<p$ that 
$$\rk H^0(C_p;Z^j)+ \rk H^0(C_p;Z^{j-1})=
\rk H^0(C_p; F_j)= {1\over p}{p\choose j}.$$
Note also that $H^0(C_p;Z^0)$ has rank 1, and so the rank of 
$H^0(C_p;Z^j)$ satisfies the same recurrence relation as $g(j)$
defined in Lemma 2.2, so is equal to $g(j)$.  Now the universal
coefficient theorem tells us that 
$$E_2^{i,j}\cong H^i(C_p;Z^j/nZ^j)\cong 
\Tor(H^{i+1}(C_p;Z^j),\Bbb Z/n)\oplus H^i(C_p;Z^j)\otimes \Bbb Z/n,$$
and we see that $E_2^{*,*}$ is as claimed.  

It will turn out that the $E_2$-page, together with the 
cohomology of $E_0^{*,*}$ (which was calculated in the last section), 
suffices to determine the pattern of higher differentials.  
We claim that the only non-zero higher differential is 
$$d_2:E_2^{i,j}\rightarrow E_2^{i+2,j-1},$$
where $i\geq 0$, $0\geq j > 1-p$ and $i+j$ is odd.  Moreover, when 
$i=0$, the kernel of $d_2$ is a direct summand of $E_2^{i,j}$.  
These claims imply those made in the statement concerning the
$E_3$-page.  Let 
$$A_m=\bigoplus_{i+j=m}E_2^{i,j},$$
so that $d_2$ and the higher differentials give rise to maps from
subquotients of $A_m$ to subquotients of $A_{m+1}$.  Then for 
$m>0$, the order of $A_m$ is $p^p$, while for $0\geq m\geq 1-p$ 
the order of $A_m$ is $n^{g(-m)}p^{2[(p-m)/2]}$ (here the square
brackets indicate the greatest integer less than or equal to 
their contents).  The order of $H^m\tot(E_0^{*,*})$ is, by Theorem~4.1, 
equal to $n^{g(-m)}p$ for $m\geq 2-p$, and $n^{g(-m)}$ for 
$m\leq 1-p$.  For $m=1-p$, the orders of $A_m$ and $H^m(E_0^{*,*})$ 
are equal, and so there are no non-zero higher differentials leaving 
$A_{1-p}$.  (This is also easy to see directly, because $A_{1-p}$ is 
concentrated in the bottom left corner of the $E_2$-page.)  Now 
assume that $m$ is even, and that we have shown that $A_m$ consists
of cycles for all the higher differentials.  In this case, the 
universal cycles in $A_{m+1}$ form a filtration of
$H^{m+1}(E_0^{*,*})$, and in particular these two groups have the
same order.  The universal cycles in $A_{m+2}$ modulo the image of 
$A_{m+1}$ under all differentials forms a filtration of
$H^{m+2}=H^{m+2}(E_0^{*,*})$.  Taking the orders of various groups we
obtain the following inequality, 
$$|H^{m+2}|= |\hbox{cycles in }A_{m+2}|/|\hbox{image of }A_{m+1}|
\leq |A_{m+2}|.|H^{m+1}|/|A_{m+1}|, $$
with equality if and only if $A_{m+2}$ consists entirely of universal
cycles.  It is easy using the numbers given above to check that 
equality does hold, and so by induction $A_{m+2}$ consists of 
universal cycles.  

Now pick any strictly positive $m$.  By the above inductive argument,
all elements of $A_{2m+2}$ are universal cycles, and it is clear that
$E_2^{2m+2,0}$ cannot be hit by any differential.  However,
$H^{2m+2}$ has order $p$, so the other $p-1$ summands of $A_{2m+2}$
of order $p$ must all be hit by some differential.  The only way for 
this to happen is if 
$$d_2:E_2^{2m+1+i,-i}\rightarrow E_2^{2m+3+i,-i-1}$$
is an isomorphism
for $0\leq i<p-1$.  The $H^*(C_p;\Bbb Z)$-module structure of the
spectral sequence allows us to deduce the claimed description of $d_2$
on $E_2^{i,j}$ for all $i>0$.  To see that $d_2$ on $E_2^{0,j}$ is
non-zero for $j$ odd and $0>j>1-p$, note that no higher differential
can hit a generator for $E_2^{2,j-1}$, but this generator must be hit
to ensure that the product over $i$ of the orders of 
$E_\infty^{2+i,j-i-1}$ equals the order of $H^{2+j-1}(E_0^{*,*})$.  
To see that the kernel of $d_2$ on $E_2^{0,j}$ is isomorphic
to $(\Bbb Z/n)^{g(-j)}$ for $j$ odd, note that since $d_2:E_2^{2,j}
\rightarrow E_2^{4,j-1}$ is an isomorphism, the kernel of
$d_2:E_2^{0,j}\rightarrow E_2^{2,j-1}$ is equal to the kernel of the 
map from $H^0(C_p;Z^j/nZ^j)=E_2^{0,j}$ to
$H^2(C_p;Z^j/nZ^j)=E_2^{2,j}$ induced by the cup product with a
generator of $H^2(C_p;\Bbb Z)$.  This kernel contains the image of 
$H^0(C_p;Z^j)$, which is isomorphic to $(\Bbb Z/n)^{g(-j)}$, because 
$H^2(C_p; Z^j)=0$ for odd $j$.  
Given the description of $E_3$ claimed in the statement, it is easy
to see that there are no further non-zero differentials, and that 
the non-trivial extension problems must be as claimed.  
\qed

\proclaim Theorem 5.1$'$.  Let $p=2$, let $E_0^{*,*}$ and
${E'}_0^{*,*}$ be as in the statement of Theorem~4.1$'$, and let 
$E_*^{*,*}$, ${E'}_*^{*,*}$ be the corresponding type I spectral
sequences.  Then if $n$ is odd, 
$$E_2^{i,j}=\cases{\relax\Bbb Z/n&if $i=j=0$,\cr
0&otherwise}$$
$${E'}_2^{i,j}=\cases{\relax\Bbb Z/n&if $i=0$, $j=1$,\cr
0&otherwise}$$
and both spectral sequences of course collapse at $E_2$.  
\pra 
If $n$ is even, then 
$$E_2^{i,j}=\cases{\relax\Bbb Z/n&if $i=j=0$,\cr
\relax\Bbb Z/2&if $j=-1$ and $i\geq 0$, or if $j=0$ and $i>0$,\cr
0&otherwise}$$
$${E'}_2^{i,j}=\cases{\relax\Bbb Z/n&if $i=0$, $j=1$,\cr
\relax\Bbb Z/2&if $j=2$ and $i\geq 0$, or if $j=1$ and $i>0$,\cr
0&otherwise.}$$
Both spectral sequences collapse at $E_3$, and 
$$E_3^{i,j}=\cases{\relax\Bbb Z/n&if $i=j=0$,\cr
\relax\Bbb Z/2&if $(i,j)=(0,-1)$ or $(1,-1)$, or if $i>0$ and even and
$j=0$ or $-1$,\cr
0&otherwise,}$$
$${E'}_3^{i,j}=\cases{\relax\Bbb Z/n&if $i=0$ and $j=1$,\cr
\relax\Bbb Z/2&if $i$ is odd and $j=1$ or 2,\cr
0&otherwise.}$$
The only non-trivial extension in reassembling $H^*\tot E_0^{*,*}$ and
$H^*\tot E_0^{\prime *,*}$ from the $E_\infty$-pages is that the
extension with kernel $E_3^{1,-1}$ and quotient $E_3^{0,0}$ is
non-split.  

\proof Similar to, but far easier than, that of Theorem~5.1.  
\qed

\beginsection 6. A DETECTION LEMMA

Let $X$ be a finite type CW-complex and 
let $\alpha$ and $\beta$ be the maps: 
$$\alpha : X^p\times EC_p\rightarrow X^p\times_{C_p}EC_p,\qquad
\beta : X\times BC_p \rightarrow X^p\times_{C_p}EC_p.$$
Here $\alpha$ is the covering map, and $\beta$ is induced by the map 
$\Delta:(x,e)\mapsto (x,\ldots,x,e)$ from 
$X\times EC_p$ to $X^p\times EC_p$,
which is $C_p$-equivariant for the trivial action on $X$ and the
permutation action on $X^p$.  One easy corollary of Nakaoka's
description of the mod-$p$ cohomology of $X^p\times_{C_p}EC_p$ is 
Quillen's detection lemma [\qui], which states that the map 
$$(\alpha^*,\beta^*): H^*(X^p\times_{C_p}EC_p;\fp) \rightarrow 
H^*(X^p\times EC_p;\fp)\times H^*(X\times BC_p;\fp)$$ 
is injective.  The result of this section is an integral analogue.  

\proclaim Corollary 6.1. Let $p$ be a prime, and as usual let $X$ be a
finite-type CW-complex, and let $\alpha$, $\beta$ be as above.  Then
the kernel of the map 
$$(\alpha^*,\beta^*): H^*(X^p\times_{C_p}EC_p) \rightarrow 
H^*(X^p\times EC_p)\times H^*(X\times BC_p)$$ 
has exponent $p$, and does not contain any cyclic summand of 
$H^*(X^p\times_{C_p}EC_p)$.  If $H^*(X)$ is expressed as a direct sum
of cyclic groups, then the kernel may be described as follows:  
\pra
For $p$ an odd prime, each cyclic summand of $H^i(X)$ of finite order
divisible by $p$ gives rise to one cyclic summand of
$\ker(\alpha^*,\beta^*)$ in each degree $pi,pi-2,\ldots,p(i-1)+3$, and
these summands form the whole of $\ker(\alpha^*,\beta^*)$.  
\pra
For $p=2$, each cyclic summand of $H^{2i}(X)$ of finite even order
gives rise to one cyclic summand of $\ker(\alpha^*,\beta^*)$ in degree
$4i$, and these summands form the whole kernel.  If $x$ generates a 
summand of $H^{2i}(X)$ of even order $n$, then $n(\wre x)$ generates 
the corresponding summand of the kernel.  

\proof Fix a splitting of $H^*(X)$ as a direct sum of cyclic groups.
Note that $\alpha^*$ is the edge map in the spectral sequence whose
$E_\infty$-page was described in Theorems 5.1~and~5.1$'$.  It follows
that $\ker(\alpha^*)$ consists of a direct sum of cyclic subgroups of
order $p$ as described in the statement, which are not direct summands
of $H^*(X^p\times_{C_p}EC_p)$, and various cyclic direct summands of
$H^*(X^p\times_{C_p}EC_p)$ of order $p$.  Consider the map of 
Cartan-Leray spectral
sequences induced by the following map of $p$-fold covering spaces:  
$$\matrix{X\times EC_p&\mapright{}& X\times BC_p\cr
\mapdown{\Delta}&&\mapdown{\beta}\cr
X^p\times EC_p&\mapright{}&X^p\times_{C_p}EC_p.\cr
}$$
The spectral sequence for $X\times BC_p$ collapses at the $E_2$-page,
and for $i>0$, $E_2^{i,j}$ consists of cyclic direct summands of
$H^*(X\times BC_p)$ of order $p$.  The map on $E_\infty$-pages is a
filtration of $\beta^*$, so it follows that the image
$\beta^*(\ker(\alpha^*))$ is a sum of cyclic summands of
$H^*(X\times BC_p)$ of order $p$.  Thus any element of $\ker(\alpha)$
not generating a cyclic summand of $H^*(X^p\times_{C_p}EC_p)$ is also
in $\ker(\beta^*)$, and $\ker(\alpha^*,\beta^*)$ is at least as large
as claimed.  To see that no cyclic summand of
$H^*(X^p\times_{C_p}EC_p)$ of order $p$ may be contained in
$\ker(\alpha^*,\beta^*)$, consider the following diagram.  
$$\matrix{
H^*(X^p\times_{C_p}EC_p)&\mapright{(\alpha^*,\beta^*)} 
&H^*(X^p)\times H^*(X\times BC_p)\cr
\mapdown{}&&\mapdown{}\cr
H^*(X^p\times_{C_p}EC_p;\fp)&\mapright{(\alpha^*,\beta^*)} 
&H^*(X^p)\times H^*(X\times BC_p;\fp)\cr
}$$
The lower horizontal map is known to be injective, and the kernel of
the left-hand map is $p.H^*(X^p\times_{C_p}EC_p)$.  Thus the lower
composite is injective on any summand of $H^*(X^p\times_{C_p}EC_p)$ of
exponent $p$, and hence so is the top horizontal map.  
\qed

\remark
One could make a more general statement for cochain
complexes $C^*$ equipped with a map from $ C^*$ to 
$C^*\otimes C^*$ having properties similar to a diagonal approximation
for the cellular cochain complex of a CW-complex.  This map is used in
the definition of the map $\beta^*$ above.

\beginsection 7. $p$-LOCAL COHOMOLOGY FOR THE SYMMETRIC GROUP 
$\Sigma_p$

The methods used in Sections 4, 5, and 6 may be used to compute the
cohomology with coefficients $\plocal$, the integers localized at $p$,
of $X^p\times_SES$ for any subgroup $S$ of the symmetric group on a 
set of size $p$.  If $S$ does not act transitively (or equivalently
has order coprime to $p$), then the Cartan-Leray type spectral sequence  
for $X^p\times_SES$ has $E_2$-page concentrated in the line
$E_2^{0,*}$.  If $S$ does act transitively, then either the results of
sections 4, 5, and 6, together with some transfer argument could be
used, or the methods used above could be applied directly.  The case
when $S$ is the full symmetric group works particularly easily.  As 
examples we give the following, which are analogues of 4.1 and
6.1.  

Throughout this section, we let $W_*$ be the chain complex with 
$\plocal$ coefficients for a 
contractible, free $\Sigma_p$-CW-complex.  Let $D'(n,i)$ be the 
complex 
$$D'(n,i)^* = D(n,i)^*\otimes \plocal = 
C(n,i)^{*\otimes p}\otimes \plocal$$
of $\plocal\Sigma_p$-modules, where $D(n,i)$ is as defined in Sections
3~and~4, and without loss of generality, $n$ may be taken to be a
power of $p$.  The complexes $D(n,2i)^{*-2pi}$ and $D(n,0)^*$ are
isomorphic, as are the complexes $D(n,2i+1)^{*-2pi}$ and $D(n,1)^*$.   

\proclaim Theorem 7.1.  Let $p$ be an odd prime, and let $n$ be a
(strictly positive) power of $p$.  Define double cochain complexes 
$E_0^{*,*}$ and $E_0^{\prime *,*}$ by 
$$E_0^{i,j}= \Hom_{\Sigma_p}(W_i,D'(n,0)^j), \qquad 
E_0^{\prime i,j}= \Hom_{\Sigma_p}(W_i,D'(n,1)^j), $$
where $W_*$ and $D'(n,i)^*$ are as above.  The corresponding type~II
spectral sequences collapse at $E_2$, and give rise to no extension
problems.  The $E_2$-pages are as follows.  
$$E_2^{i,j}= \cases{
\plocal/(pn)&for $i=j=0$,\cr
\plocal/(p)&for $i=0$, $j=(2j'+2)(p-1)$,\cr
          &\quad or $i=-p$, $j=(2j'+1)(p-1)$,\cr 
0&otherwise.\cr}  
$$
$$E_2^{\prime i,j}=\cases{\plocal/(n)&for $i=0$, $j=1$,\cr 
\plocal/(p)&for $ i=0$, $j=(2j'+2)(p-1)$,\cr
          &\quad or $i=p$, $j=(2j'+1)(p-1)$,\cr 
0&otherwise.\cr}
$$

\proof Similar to the proof of Theorem~4.1.  To calculate the
$E_1$-pages of the spectral sequences, note that as a
$\plocal\Sigma_p$-module, 
$$D(n,0)^{-j} \cong D(n,1)^j\cong
\ind^{\Sigma_p}_{\Sigma_j\times\Sigma_{p-j}}(\epsilon_j),$$ 
where $\epsilon_j$ is the sign representation for $\Sigma_j$ tensored
with the trivial representation of $\Sigma_{p-j}$.  Using the
Eckmann-Shapiro lemma it follows that 
$$E_1^{0,0}\cong E_1^{0,-1}\cong E_1^{\prime 0,0}\cong 
E_1^{\prime 0,1}\cong \plocal,$$ 
and that with these exceptions $E_1^{i,j}= E_2^{i,j}$ and 
$E_1^{\prime i,j}= E_2^{\prime i,j}$ as described 
in the statement.  As in Theorem~4.1, the collapse at $E_2$ follows by
comparing with the case $n=1$ and applying Lemma~1.1.  
\qed

\proclaim Corollary 7.2.  Let $X$ be a CW-complex of finite type, and
let $\alpha$, $\beta$ be the maps 
$$\eqalign{ X^p\times E\Sigma_p&\mapright{\alpha}
X^p\times_{\Sigma_p}E\Sigma_p,\cr
X\times B\Sigma_p&\mapright{\beta}
X^p\times_{\Sigma_p}E\Sigma_p,\cr}$$
analogous to the maps occurring in the statement of Corollary~6.1.  If
$\alpha^*$, $\beta^*$ are the induced maps on cohomology with
$\plocal$ coefficients, then $\ker(\alpha^*,\beta^*)$ is concentrated
in degrees divisible by $2p$, and may be described as follows.  Every
cyclic summand of $H^{2i}(X;\plocal)$ of order $p^r$ gives rise to a
summand of $\ker(\alpha^*,\beta^*)$ of order $p$ in degree $2pi$.  If
$x$ generates a summand of $H^{2i}(X;\plocal)$ of order $p^r$, then 
$p^r(\wre x)$ generates the corresponding summand of the kernel.  

\proof Similar to that of Corollary~6.1. \qed

\remark Another generalization of the material in Section~4 is to
consider other possibilities for $W_*$.  We have made calculations
similar to those in Sections 4~and~5 for the case when $S$ is cyclic
of order $p$, and $W_*$ is the
chain complex of a sphere with a free $S$-action.  In this case the
non-zero entries on the $E_2$-page of the type~II spectral sequence
are concentrated at the edges of a rectangle (of height $p$ and width 
the dimension of the sphere), and the left-hand edge of the rectangle
is identical to the start of the \lq infinite rectangle' considered in
Section~4.  This gives enough information to reconstruct the whole 
$E_2$-page and the higher differentials may be calculated using
Lemma~1.1.

\beginsection 8. THE EXPONENT OF GROUP COHOMOLOGY

In this section we apply the results of previous sections to the
question of determining the exponent of the integral cohomology of 
finite groups.  It is simpler to concentrate on $p$-groups, and 
we shall do this, although we shall also make some remarks about the
case of arbitrary finite groups.  

Before starting, we recall the Evens-Venkov theorem ([\evb,\ben]),
which states that if $G$ is a finite group and $H$ is a subgroup of
$G$, then $H^*(G)$ is a finitely generated ring, and $H^*(H)$ is a
finitely generated module for $H^*(G)$.  

\proclaim Definition A.  For $G$ a finite group, say that $G$ enjoys
property~A when for all $i$, if there exists $j$ such that $H^j(G)$
contains an element of order $p^i$, then there exist infinitely many
such $j$.  

Property~A was introduced by A. Adem in [\ade], and by H.-W. Henn in
[\hwh], but see also [\oprob, q.\ no.~754].  No 
$p$-group is known which does not have property~A, and Adem made 
conjecture~A; the conjecture that all $p$-groups have property~A{}.  
Henn also asked whether all $p$-groups have 
property~A{}.  If $G$ is a finite group with Sylow subgroup $G_p$, 
then $H^j(G)$ contains elements of order $p^i$ for
infinitely many $j$ if and only if $H^j(G_p)$ does.  Thus $G$ enjoys
propery~A if its Sylow-$p$ subgroup does.  

One may reformulate property~A in terms of the cohomological exponent 
$\ce(G)$ and eventual cohomological exponent $\ece(G)$ of a finite group
$G$, which we define as follows, where $\exp(-)$ stands for the
exponent of an abelian group:  
$$\ce(G) = \exp\bigl(\bigoplus_{i>0}H^i(G)\bigr),\qquad
\ece(G) = \lim_{j\rightarrow \infty}
\exp\bigl(\bigoplus_{i>j}H^i(G)\bigr).$$
Note that for any $G$, $\ece(G)$ divides $\ce(G)$ and $\ce(G)$ divides
$|G|$.  The group $G$ enjoys property~A if and only if $\ce(G) = \ece(G)$.
Adem did not make the above definition, but he pointed out that the
Evens-Venkov theorem implies that if $H\leq G$, then $\ece(H)$ divides
$\ece(G)$ [\ade].  It seems to be unknown whether a similar property
holds for $\ce(G)$, so we make 

\proclaim Conjecture A$^-$.  For $G$ any $p$-group and $H$ any subgroup of
$G$, $\ce(H)$ divides $\ce(G)$.  

We call this conjecture~A$^-$ because it is weaker than conjecture~A{}.

\proclaim Proposiion 8.1.  Let $G$ be a $p$-group.  Then 
$$\ece(G\wreath C_p) = p\cdot \ece(G)\quad \hbox{and}\quad 
\ce(G\wreath C_p) = p\cdot \ce(G),$$
except that possibly $\ce(G\wreath C_2)= \ce(G)$ if $p=2$ and $H^*(G)$
contains only finitely many cyclic summands of order $\ce(G)$, all of
which occur in odd degree.  

\proof First consider the case when $p$ is odd.  Express $H^*(G)$ as a
direct sum of cyclic subgroups.  All of these have order dividing
$|G|$, except that $H^0(G) =\Bbb Z$.  By Theorem~4.2, each $p$-tuple
of cyclic summands of $H^*(G)$, not all equal, of orders
$n_1,\ldots,n_p$ gives rise to finitely many cyclic summands of
$H^*(G\wreath C_p)$ of order ${\rm h.c.f.}\{n_1,\ldots,n_p\}$.  Each cyclic
summand of $H^*(G)$ of order $p^i$ also contributes $(p-1)/2$ cyclic
summands of $H^*(G\wreath C_p)$ of order $p^{i+1}$, infinitely many
summands of order $p$, and finitely many other summands of order
dividing $p^i$.  The only other contribution is from $H^0(G)$, which
gives rise to $H^0(G\wreath C_p)$, and other summands of order $p$.  
The claim follows.  

The case when $p=2$ is similar, relying on Theorem~4.2$'$.  The extra
difficulty arises because summands of $H^*(G)$ of order $2^i$ in odd
degree do not contribute any summands of $H^*(G\wreath C_2)$ of order
$2^{i+1}$.  It remains to rule out the possibility that $H^*(G)$
contains infinitely many cyclic summands of order $2^i=\ece(G)$, all but
finitely many of which occur in odd degrees.  However, $H^*(G)$ is a
finitely generated ring, and since the square of any generator has
even degree, $H^*(G)$ is a finitely generated module for the subring
of elements of even degree.  Let the supremum of the degrees of a
(finite) set of module generators for $H^*(G)$ over the even degree
subring be $m$.  Now suppose that all elements of $H^{2j}(G)$ have 
exponent strictly less than $2^i$ whenever $j> j_0$.  In this case, 
$H^j(G)$ may contain elements of exponent $2^i$ only for $j\leq
2j_0+m$, and so $H^*(G)$ can have only finitely many cyclic summands
of order $2^i$.  \qed

\proclaim Corollary 8.2.  Let $G$ be a $p$-group.  If $p$ is odd, 
then $G$ enjoys property~A if and only if $G\wreath C_p$ does.  
If $p=2$ and $G$ enjoys property~A, then so does $G\wreath C_2$.
\qed

The following slightly stronger property than property~A fits well
with wreath product arguments, as can be seen from Theorem~8.3.  

\proclaim Definition A$'$.  For $p$ an odd prime, say that a $p$-group
$G$ enjoys property~A$'$ if for each $i>0$, whenever there exists $j$ such that
$H^j(G)$ contains a cyclic summand of order $p^i$, there exist
infinitely many such $j$.  
\pra
Say that a 2-group $G$ enjoys property~A$'$ if for each $i>0$ and $\epsilon
= 0,1$, whenever there exists a $j$ such that $H^{2j+\epsilon}(G)$
contains a cyclic summand of order $p^i$, there exist infinitely many
such $j$.  

I know of no $p$-groups not enjoying property~A$'$, and it seems
reasonable to make conjecture~A$'$, i.e., to conjecture that all
$p$-groups enjoy property~A$'$.  

\proclaim Theorem 8.3.  Let $G$ and $S'$ be $p$-groups, and let $S'$
act on the finite set $\Omega$, with image $S\leq \Sigma(\Omega)$.  If
both $G$ and $S'$ enjoy property~A$'$, then so does $G\wreath S'$.  

\proof As at the end of Section~1, we fix a splitting of $H^*(G)$ as a
direct sum of cyclic groups indexed by a set $A$, and let 
$$C'^*=\bigoplus_{a\in A} C(a)^*$$ 
be a direct sum of complexes of the form $C(n,i)^*$ having cohomology
isomorphic to $H^*(G)$.  

Let $\bf a$ stand for the $S'$-orbit (or equivalently $S$-orbit) in 
$A^\Omega$ containing $(a_1,\ldots,a_l)$, and let $D({\bf a})^*$ be
the $\Bbb ZS$-subcomplex of $C'^{*\otimes\Omega}$ generated by 
$C(a_1)\otimes\cdots\otimes C(a_l)$.  Then as complexes of $\Bbb
ZS$-modules, 
$$C'^{*\otimes \Omega}= \bigoplus_{{\bf a}\in A^\Omega/S} D({\bf
a})^*,\eqno{(3)}$$
where the sum is over the $S$-orbits in $A^\Omega$.  
For $p$ odd, the isomorphism type of $D({\bf a})^*$ (modulo a shift in
degree) depends only on the orders of
the summands $H(a_1),\ldots,H(a_l)$ of $H^*(G)$, because no subgroup
of $S'$ has a non-trivial sign representation.  For $p=2$, the
isomorphism type of $D({\bf a})^*$ (modulo an even shift in degree) 
depends on the parity of the
degrees of the $H(a_i)$ as well as on their orders.  In either case,
property~A$'$ for $G$ implies that each isomorphism type of $D({\bf
a})^*$ that occurs in $(3)$ occurs infinitely often, except for 
$D({\bf a}_0)$, ${\bf a}_0 = (a_0,\ldots,a_0)$, where $H(a_0) =
H^0(G)\cong \Bbb Z$.  The complex $D({\bf a}_0)$ is isomorphic to the
trivial $S$-module $\Bbb Z$ concentrated in degree zero.  

Let $K$ be the kernel of the homomorphism from $S'$ onto $S$, and let
$W_*$ be the cellular chain complex for $ES'/K$, viewed as a complex
of free $\Bbb ZS$-modules.  By Lemma~1.4, $H^*(G\wreath S')$ is
isomorphic to the cohomology of the total complex of 
$$\eqalign{E_0^{*,*}&=\Hom_S(W_*,\bigoplus_{{\bf a}\in
A^\Omega/S}D({\bf a})^*)\cr
&=\bigoplus_{{\bf a}\in A^\Omega/S}\Hom_S(W_*,D({\bf a})^*)\cr
&=\bigoplus_{{\bf a}\in A^\Omega/S}E({\bf a})_0^{*,*}.\cr}$$
\par From the analysis of the $D({\bf a})^*$'s given above, it follows that
each isomorphism type of $E({\bf a})_0^{*,*}$ that occurs, occurs
infinitely often, except for $E({\bf a}_0)_0^{*,*}$.  But 
$$E({\bf a}_0)_0^{*,*}\cong \Hom_S(W_*,\Bbb Z),$$
and so $H^*\tot E({\bf a}_0)_0^{*,*}$ is isomorphic to $H^*(S')$.
Thus if ${\bf a}\neq {\bf a}_0$, then any cyclic summand of 
$H^*\tot E({\bf a})_0^{*,*}$ is also a summand of 
$H^*\tot E({\bf a}')_0^{*,*}$ for infinitely many ${\bf a}'$ by
property~A$'$ for $G$, while any cyclic summand of $H^*\tot E({\bf
a}_0)_0^{*,*}$ occurs infinitely often by property~A$'$ for $S'$.  
\qed

For $G$ a $p$-group, define $\Phi_n(G)$ to be the intersection of the
subgroups of $G$ of index at most $p^n$.  Thus $\Phi_1(G)$ is the
Frattini subgroup of $G$, and each $\Phi_n(G)$ is a characteristic
subgroup of $G$.  In [\ian] we pointed out that these subgroups give a
group-theoretic description of a good upper bound for $\ece(G)$, as
follows.  

\proclaim Proposition 8.4.  Let $G$ be a $p$-group such that
$\Phi_n(G)=\{1\}$.  Then $\ece(G)$ divides $p^n$.  

\proof If $H_1,\ldots,H_m$ are a family of subgroups of $G$ with
trivial intersection, then the natural map 
$$G\longrightarrow \Sigma(G/H_1)\times\cdots\times\Sigma(G/H_m)$$
sending an element $g$ to the permutation $g'H_i \mapsto gg'H_i$ 
is injective.  If $G$ is a $p$-group such that $\Phi_n(G)$ is trivial,
then $G$ is therefore isomorphic to a subgroup of a product of copies
of the Sylow $p$-subgroup $P_n$ of $\Sigma_{p^n}$.  Using the results
of Section~4, it may be shown that $\ece(P_n) = \ce(P_n) = p^n$, and the 
claim follows.  \qed

\remarks It may be shown that $\ece(P_n) = \ce(P_n) = p^n$ without using
the results of Section~4---the upper bound by a transfer argument,
and the lower bound by exhibiting a suitable subgroup of $P_n$ whose
cohomology is known, for example the cyclic group of order $p^n$.  

If $G$ is a finite group with Sylow $p$-subgroup $G_p$, then the
$p$-part of $\ece(G)$ is equal to $\ece(G_p)$, so the bound given by 
Proposition~8.4 may be applied to arbitrary finite groups.  

The bound on $\ece(G)$ given by Proposition~8.4 is sharp in many cases
including:  the extraspecial groups, and the groups
presented during the proof of Theorem~3.3; groups of $p$-rank one; 
abelian groups; the metacyclic groups whose cohomology has been
calculated by Huebschmann in [\huea].  Moreover, if $G$ 
has subgroups $H_1,\ldots,H_m$ of
index $p^i$ with trivial intersection then the subgroups of the form 
$G\times\cdots\times H_i\times \cdots\times G$ of $G^p$
intersect trivially and have index $p^{i+1}$ in $G\wreath C_p$.  
It follows from Proposition~8.1 that whenever the bound on $\ece(G)$ is sharp,
so is the bound on $\ece(G\wreath C_p)$.  In the next section we shall 
exhibit a group for which the bound is not sharp however.  This is a 
2-group $G$ such that $\ece(G)=4$ but $\Phi_2(G)\neq\{1\}$.

\beginsection 9. A GROUP WHOSE COHOMOLOGY HAS SMALL EXPONENT

As promised at the end of Section~8, we exhibit a counter-example to
the converse of Proposition~8.4.  More precisely, we exhibit a 2-group 
$G$ whose index four subgroups do not intersect trivially, but such 
that $\ece(G)=4$.  It is known that if $H$ is a 
$p$-group such that $\ece(H)=p$,
then $H$ is elementary abelian, and so $\Phi_1(H)$ is trivial
[\ade,\hwh,\ian].  Thus our example is minimal in some sense.  
The smallest such $G$ that we have been able to find has order $2^7$.
The necessary calculations are simpler for another example of order 
$2^9$ however, so we shall concentrate on this example and explain the
smaller example in some final remarks.  We also discuss the case of 
odd $p$ at the end of the section.  

To explain the origin of the example, it is helpful to
consider a more general problem.  If $\Gamma$ is a discrete group of
finite virtual cohomological dimension (vcd), then $\ece(\Gamma)$ may be
defined just as for finite groups.  Just as in Proposition~8.4, it may
be shown that if the subgroups of $\Gamma$ of index dividing $p^i$
have torsion-free intersection, then $\ece(\Gamma)$ divides $p^i$.  In
this case, groups $\Gamma$ for which this bound is not tight are
already known.  Let $\Gamma(n)$ (depending on the prime $p$ as well as
the integer $n$) be the group with presentation:  
$$\eqalign{\Gamma(n)=\langle A_1,\ldots,A_n,B_1,\ldots,B_n,C \mid
C^p,\, &[A_i,C],\, [B_i,C],\, \cr
&[A_i,A_j],\, [B_i,B_j],\,
[A_i,B_j]C^{-\delta(i,j)}\rangle.\cr}$$
Thus $\Gamma(n)$ is expressible as a central extension with kernel 
cyclic of order $p$ generated by $C$ and quotient free abelian of rank
$2n$.  The spectral sequence with $\Bbb Z$-coefficients for this 
central extension collapses.  By examining this spectral sequence 
Adem and Carlson were able to determine the ring 
$H^*(\Gamma(n))$ [\adca].  From their calculation it follows 
that $\ece(\Gamma(p))= p^p$.  

On the
other hand, it is easy to show that $\Phi_p(\Gamma)$, the intersection
of the subgroups of $\Gamma$ of index dividing $p^p$, contains the
element $C$, so is not torsion-free.  One way to see this is to note 
that the centre $Z(\Gamma)$ is generated by $C$, $A_i^p$, and $B_i^p$,
while the commutator subgroup $\Gamma'$ of 
$\Gamma$ is generated by $C$.  The quotient $\Gamma/Z(\Gamma)$ 
is an elementary abelian $p$-group of rank $2p$.  The map 
$$\Gamma/Z(\Gamma)\times \Gamma/Z(\Gamma)\rightarrow \Gamma'$$
given by $(g,h)\mapsto [g,h]$ may be viewed as an alternating bilinear
form on the $\fp$-vector space $\Gamma/Z(\Gamma)$, whose maximal
isotropic subspaces have dimension $p$.  Now any subgroup of $\Gamma$
of index at most $p^p$ either contains $Z(\Gamma)$, and in particular 
contains $C$, or has image in $\Gamma/Z(\Gamma)$ a
subspace of dimension greater than $p$, which cannot be an isotropic
subspace.  Hence $C$ is a commutator in any subgroup of index $p^p$ that 
does not contain $Z(\Gamma)$.

The quotient of $\Gamma$ by the subgroup generated by the $A_i^p$'s
and $B_i^p$'s is extraspecial of order $p^{2p+1}$ (in the case $p>2$,
this group has exponent $p$, and in the case $p=2$ it is a central
product of copies of the dihedral group of order eight).  For this
group it is known that the bound on the eventual cohomological
exponent given by Proposition~8.4 is best possible (see for example 
Lemma~3.5).  Our example is a
slightly larger quotient of $\Gamma(2)$ in the case when $p=2$.  
Let $G$ be the group with presentation 
$$
\eqalign{G=\langle A_1,A_2,B_1,B_2,C \mid C^2,\, &A_i^4,\, B_i^4,\, 
[A_i,C],\cr
\, &[B_i,C],\, [A_1,A_2],\, [B_1,B_2], [A_i,B_j]C^{\delta(i,j)}
\rangle.}\eqno{(4)}$$
Thus $G$ is the quotient of $\Gamma(2)$ by the subgroup generated by 
$A_i^4$ and $B_i^4$.  The subgroup $Z$ of $G$ generated by the $A_i^2$ and
$B_i^2$ is central, and elementary abelian of rank four.  The quotient
$Q=G/Z$ is the extraspecial group of order $2^5$ consisting of a
central product of two copies of the dihedral group of order eight.  

Our partial calculation of $H^*(G)$ relies on some knowledge of the
ring structure of $H^*(Q)$, at least up to degree four.  The additive 
structure of $H^*(Q)$, and of the integral cohomology of all 
extraspecial 2-groups, was determined by Harada and Kono 
[\hako,\bencarl].  We find that the description of the 
ring structure of $H^*(Q)$ given in [\bencarl] is incorrect, 
even in low degrees.  All that we
shall require concerning $H^*(Q)$ is contained in the following
statement.  

\proclaim Lemma 9.1.  Let $Q$ be the quotient $G/Z$ as above, where 
$G$ has the presentation (4).  There is an element $\chi$ of $H^4(Q)$ of
order eight, and $4\chi$ is expressible as a sum of products of
elements of $H^2(Q)$.  
%\pra 
%More precisely, let $y_1,\ldots,y_4$ form a basis for $H^1(E;\ftwo)$,
%and let $y_1y_2+y_3y_4=0$ in $H^2(E;\ftwo)$.  Let $\delta$ stand for
%the Bockstein from mod-2 cohomology to integral cohomology, and let 
%$x_i = \delta(y_i)$.  Then $4\chi= x_1x_2+x_3x_4$.   

\proof First we show that there is an element of $H^4(Q)$ of order
eight (we could also quote this fact from [\hako]).  Recall that 
$Q$ has seventeen irreducible real representations, sixteen 
1-dimensional ones and one 4-dimensional.  Let $S$ be the
unit sphere in the 4-dimensional faithful real representation.  
Then $Q$ acts trivially on $H^*(S)$, so there is an Euler class
$e(S)\in H^4(Q)$ defined for $S$.  (Topologically, $e(S)$ is the Euler
class of the orientable $S$-bundle $S\times_QEQ$ over $BQ= EQ/Q$.
Algebraically, $E(S)$ is the extension class in 
$\Ext^4_{\Bbb ZQ}(\Bbb Z,\Bbb Z)$ represented by the chain complex of
the $Q$-CW-complex $S$.)  The centre $Z(Q)$ of $Q$ is cyclic of order two 
generated by the image of $C$, and the 4-dimensional faithful real
representation restricts to this subgroup as four copies of its
non-trivial real representation.  It follows that the image of $e(S)$
generates $H^4(Z(Q))$.  Now $Q$ is a group of the type
discussed in Lemma~3.5, and this lemma implies that $e(S)$ has order
eight.  

For the rest of the proof, we consider the spectral sequences for the 
central extension 
$$Z(Q) \rightarrow Q \rightarrow Q/Z(Q)\eqno{(5)}$$
with integer and mod-2 coefficients.  The spectral sequence with mod-2
coefficients was solved completely by Quillen [\quies,\bencarl].  
Recall that the mod-2 cohomology ring of an elementary abelian 2-group
of rank~$r$ is a polynomial ring on $r$~generators of degree one.  
Let $E_*^{\prime *,*}$ be the spectral sequence for (5) with mod-2
coefficients.  Thus $E_2^{\prime *,*}= \ftwo[y_1,\ldots,y_4,z]$, where
$y_i\in E_2^{\prime 1,0}$ and $z\in E_2^{\prime 0,1}$.  The extension
class for $Q$ in $H^2(Q/Z(Q);\ftwo)$ is (without loss of generality) 
$y_1y_2+y_3y_4$.  Thus $d_2(z) = y_1y_2+y_3y_4$.  We shall not need to
consider any higher differentials in this spectral sequence.  

Now let $E_*^{*,*}$ be the spectral sequence for (5) with $\Bbb
Z$-coefficients.  Then $E_2^{i,j}= H^i(Q/Z(Q);H^j(Z(Q)))$ is more
complicated, because the integral cohomology of $Q/Z(Q)$ is more
complicated.  Except in degree zero, $H^*(Q/Z(Q))$ has exponent two,
so the Bockstein map $\delta$ and projection map $\pi$
$$\delta : H^*(Q/Z(Q);\ftwo)\rightarrow H^{*+1}(Q/Z(Q))
\qquad \pi : H^*(Q/Z(Q))\rightarrow H^*(Q/Z(Q);\ftwo)$$ 
are (respectively) surjective in positive degrees, and injective in
positive degrees.  As a ring, $H^*(Q/Z(Q))$ is generated by four
elements $x_i$ of degree two, six elements $w_{ij}$ of degree three, 
four elements $v_{ijk}$ of degree four, and one element $u=u_{1234}$ 
of degree five.  Here the indices satisfy $1\leq i<j<k\leq 4$, and 
in terms of $\delta$ the elements are
$$x_i = \delta(y_i),\qquad w_{ij}= \delta(y_iy_j),\qquad 
v_{ijk} = \delta(y_iy_jy_k), \qquad u = \delta(y_1y_2y_3y_4).$$

Note that $E_2^{1,3}$ and $E_2^{3,1}$ are trivial, and that except for
$E_2^{0,0}$, each group on the $E_2$-page has exponent two.  It
follows that if $\chi\in H^4(Q)$ has order eight, then $4\chi$ 
yields an element of $E_\infty^{4,0}$ in the spectral sequence.  The
proof will be complete once we have shown that $E_\infty^{4,0}$ is
generated by products of elements of $E_\infty^{2,0}$.  Let $t$ be a
generator for $E_2^{0,2}$.  It is logical to denote the generators of
$E_2^{1,2}$ by $[ty_1],\ldots,[ty_4]$.  Although there is no element 
$y_i$ in $E_2^{1,0}$, such relations as $[ty_1]^2=t^2x_1$ do hold in
$E_2^{*,*}$.  The group $E_2^{4,0}$ is elementary abelian of rank
fourteen, with generators the ten monomials in the $x_i$'s, and the
four elements $z_{ijk}$.  The only differential that can hit this 
group is $d_3: E_3^{1,2} \rightarrow E_3^{4,0}$.  To compute this
differential, we recall that there is a Bockstein map of spectral
sequences:  
$$\delta : E_r^{\prime i,j}\rightarrow 
\cases{ E_{r+1}^{i,j+1}&for $j>0$, \cr
E_{r+1}^{i+1,j}&for $j=0$,\cr}$$ 
such that the induced map on $E_\infty$-pages is a filtration
of the Bockstein $\delta : H^*(Q;\ftwo) \rightarrow H^{*+1}(Q)$.  
The map on $E_2^{i,j}$ is given by the map $\delta :
H^j(Z(Q);\ftwo)\rightarrow H^{j+1}(Z(Q))$ for $j>0$ and by 
$\delta : H^i(Q/Z(Q);\ftwo) \rightarrow H^{i+1}(Q/Z(Q))$ for $j=0$.  

Now the differential $d_2 : E_2^{\prime 1,1}\rightarrow E_2^{\prime
3,0}$ satisfies $d_2(zy_i)= (y_1y_2+y_3y_4)y_i$,  and so in
$E_*^{*,*}$, 
$$\eqalign{
d_3([ty_i]) = d_3\delta(zy_i) &= \delta ((y_1y_2+ y_3y_4)y_i) \cr
&= \cases{x_1x_2+z_{134}&for $i=1$,\cr
x_1x_2+z_{234}&for $i=2$,\cr
x_3x_4+z_{123}&for $i=3$,\cr
x_3x_4+z_{124}&for $i=4$.\cr}\cr}
$$
It follows that $E_4^{4,0}$ is generated by monomials in the $x_i$'s,
as required.  
\qed

\proclaim Corollary 9.2.  The group $G$ with presentation (4) above 
has $\Phi_2(G)\neq \{1\}$ and $\ece(G)=4$.  

\proof The proof that $\Phi_2(G)$ contains the element $C$ (and is in
fact equal to the subgroup generated by $C$) is identical to the proof
given above for the group $\Gamma(n)$, so shall be omitted.  

The commutator subgroup $G'$ of $G$ is cyclic of order two generated by
$C$, and $G/G'$ is isomorphic to a product of four
cyclic groups of order four.  Thus $G/G'$  has four 1-dimensional
complex representations whose kernels intersect in the trivial group,
and hence a free action with trivial action on homology on the torus 
$U(1)^4$.  Letting $S$ be the 3-sphere as in the proof of Lemma~9.1,
with $G$ acting via the action of its quotient $Q=G/Z$, it follows
that $G$ acts freely with trivial action on homology on $S\times
U(1)^4$.  By Venkov's proof of the Evens-Venkov theorem [\ben], it 
follows that $H^*(G)$ is finite over the subring generated by the 
Euler classes for these five $G$-spheres.  Thus it suffices to show
that the orders of these Euler classes are (at most) four.  

$H^2(G)\cong H^2(G/G')$ is isomorphic to 
four copies of $\Bbb Z/4$, and the Euler
classes of the four $U(1)$'s generate $H^2(G)$.  Thus it remains to
check that $e(S)$ has order four in  
$H^*(G)$.  From Lemma~9.1 we know that $e(S)\in H^*(G/Z) $ has order 
eight, but $4 e(S)$ is a sum of products of elements of
$H^2(G/Z)=(\Bbb Z/2)^4$.  The image of $H^2(Q)$ in $H^2(G)$ consists
of the elements $2x$ for $x\in H^2(G)= (\Bbb Z/4)^4$, and a product
of any two such elements $(2x)(2x')= 4xx'$ is zero in $H^4(G)$.  \qed

\remarks With a little more work it can be shown that, with notation
as in the proof of Lemma~9.1, $4\chi = x_1x_2+x_3x_4$.  Let 
$G_2 = G/\langle B_1^2,B_2^2\rangle$, so that $G_2$ is a quotient of
$G$ of order $2^7$, and $Q$ is a quotient of $G_2$.  
The argument of Corollary~9.2 can be used to show that the image of
$4\chi$ in $H^4(G_2)$ is zero, and so $\ece(G_2)= 4$.  

The evidence of Adem-Carlson's infinite groups suggests that there
should be similar examples for all primes of groups with $\ece(G)=p^p$
and $\Phi_p(G)\neq \{1\}$.  The smallest candidates (of the same type
as our example for $p=2$) have order
$p^{3p+1}$, and I have been unable to compute $\ece(G)$ for any of
these groups for $p$ odd.

\beginsection 10. VARIETIES FOR HIGHER TORSION IN COHOMOLOGY 

The description of $H^*(X^p\times_{C_p}EC_p)$ given in Sections 4--6
is complicated.  In this section we shall extract some information 
of a more conceptual nature.  We assume that $H^*(X)$ is a finitely
generated ring, and describe the ring $H^*(X^p\times_{C_p}EC_p)$ from
the point of view of algebraic geometry.  We start by defining the 
varieties that we shall study, and stating some of their elementary 
properties, before stating our main result as Theorem~10.6.   
It is convenient to use cohomology with
$\plocal$-coefficients, so throughout this section the coefficients
for cohomology are $\plocal$ when omitted.  

Fix an algebraically closed field $k$ of characteristic $p>0$.  Let
$R$ be a finitely generated commutative $\plocal$-algebra, and define
$V(R)$ to be the variety of all ring homomorphisms from $R$ to $k$,
with the Zariski topology.  For each $j\geq 0$, let $I_j= I_j(R)$ be
the annihilator in $R$ of the element $p^j$, or equivalently 
the ideal of $R$ generated by the elements of order dividing $p^j$: 
$$I_j = \sum_{p^jr=0} Rr = \Ann_R(p^j). $$
Let $V_j(R)$ be the corresponding subvariety of $V(R)$.  Since $R$
is Noetherian, there exists $j_0$ such that the $I_j$ are all equal
for $j\geq j_0$, and we define $V_\infty(R) = V_{j_0}(R)$.  Thus the 
$V_j(R)$'s are included in each other as follows:  
$$V(R)= V_0(R)\supseteq V_1(R) \supseteq \cdots \supseteq V_{j_0}(R)
= V_\infty(R).$$
A ring homomorphism $f: R\rightarrow S$ induces maps $f^*: V_j(S)
\rightarrow V_j(R)$ for all $j$.  

For any finite type CW-complex~$X$, 
let $H^\spot(X)$ stand for the subring of $H^*(X)$ consisting of
elements of even degree.  If $H^*(X)$ is a finitely generated 
$\plocal$-algebra, then $H^\spot(X)$ is an algebra of the type
considered above.  The main result of the section,
Theorem~10.6, gives a description of
$V_i(H^\spot(X^p\times_{C_p}EC_p))$ for such $X$ in terms of the 
$V_j(H^\spot(X))$'s.  Note that the case $i=0$ is not a special case
of Quillen's work on equivariant cohomology rings [\quib], because 
we do not assume that $X$ is finite.  The case $i=0$ does however
follow from Nakaoka's description of $H^*(X^p\times_{C_p}EC_p;\fp)$
together with some of the commutative algebra from the appendix to
[\quib].  The idea of studying the $V_i$ for $i>0$ came from a paper
of Carlson, who considered a special case which we shall discuss in the
next section [\carl].  Proposition~10.1 is Proposition~B.8 of [\quib].

\proclaim Proposition 10.1. ([\quib]) Let $f:R\rightarrow S$ be a
homomorphism of finitely generated commutative $\fp$-algebras, and
assume that $f$ has the following properties:  
\pra 
i) For $r\in\ker(f)$ there exists $n$ such that $r^n=0$; 
\pra 
ii) For $s\in S$, there exists $n$ such that $s^{p^n}\in \im(f)$.  
\pra\noindent
Then $f$ induces a homeomorphism from $V(S)$ to $V(R)$.  
\qed

\remark A homomorphism $f$ having properties i) and ii) is known as an
F-isomorphism.  

\proclaim Proposition 10.2. Let $X$ be a finite-type CW-complex.  Then
the following are equivalent.
\pra
\item{i)} $H^*(X)$ is a finitely generated $\plocal$-algebra; 
\pra 
\item{ii)} $H^*(X;\fp)$ is a finitely generated $\fp$-algebra, and   
the torsion in $H^*(X)$ has bounded exponent.  

\proof Let $I_j$ be the ideal of $H^*(X)$ of elements annihilated by
$p^j$.  Under hypothesis i), $H^*(X)$ is Noetherian, so there exists 
$j_0$ such that $I_j=I_{j_0}$ for all $j\geq j_0$.  Now the map from 
$H^*(X)$ to $H^*(X;\fp)$ has image $H^*(X)/(p)$, a finitely generated 
$\fp$-algebra, and cokernel isomorphic to the ideal $I_1$, which is a
finitely generated module for $H^*(X)$, and hence also for
$H^*(X)/(p)$.  Thus i) $\Rightarrow$ ii).  

Conversely, assume that ii) holds.  Let $p^j$ annihilate the torsion
in $H^*(X)$, and let $R_j$ be the image of the map $x\mapsto
x^{p^{j+1}}$ from $H^*(X;\fp)$ to itself.  Then $R_j$ is a finitely
generated ring, and $H^*(X;\fp)$ is a finitely generated $R_j$-module.
Since $p^j$ annihilates torsion in $H^*(X)$, the Bockstein spectral
sequence with $E_1$-page $E_1^i=H^i(X;\fp)$ collapses at the
$E_{j+1}$-page, so that any element of $H^*(X;\fp)$ which is a cycle
for the Bockstein $\beta_1$ and for each higher Bockstein
$\beta_2,\ldots,\beta_j$ is in the image of $H^*(X)$.  But if $x\in
H^*(X;\fp)$ is a cycle for $\beta_i$, then 
$$\beta_{i+1}(x^p)= p\beta_{i+1}(x)x^{p-1}=0,$$
so by induction, the subring $R_j$ consists of universal cycles in the
Bockstein spectral sequence.  

Moreover, each of the non-zero differentials in the Bockstein spectral
sequence is an $R_j$-linear map.  Hence the universal cycles in the
Bockstein spectral sequence, or equivalently, the image of $H^*(X)$ in
$H^*(X;\fp)$, form an $R_j$-submodule of $H^*(X;\fp)$.  It follows
that this image, which is isomorphic to $H^*(X)/(p)$, is a finitely
generated ring.  Now take a finite set of elements of $H^*(X)$ mapping
to a set of generators for $H^*(X)/(p)$.  The $\plocal$-subalgebra of
$H^*(X)$ generated by these elements is the whole of $H^*(X)$, since a
proper $\plocal$-submodule of the finitely generated module $H^i(X)$
cannot have image the whole of $H^i(X)/(p)$.  
\qed

\remark The case when $X$ is a Moore space with $H_n(X;\Bbb Z)$
isomorphic to the rationals shows that the assumption that $X$ has
finite type is necessary in Proposition~10.2.  I do not know of any 
finite type $X$ such that $H^*(X;\fp)$ is finitely generated,
but $H^*(X)$ is not finitely generated.  

\proclaim Proposition 10.3.  Let $X$ and $X'$ be CW-complexes of
finite type whose $\plocal$-cohomology is a finitely generated
algebra, and let $i\geq 0$ be a positive integer.  Then the obvious maps
induce homeomorphisms as shown below:  
\pra 
i) $V_0(H^\spot(X;\fp))\cong V_0(H^\spot(X))$;
\pra
ii) if $p=2$, then $V_0(H^*(X;\ftwo))\cong V_0(H^\spot(X;\ftwo))$; 
\pra 
iii) $V_i(H^\spot(X\times X'))\cong V_i(H^\spot(X)\otimes
H^\spot(X'))$.

\proof Firstly, note that for any $R$, the natural
map $R\rightarrow R/(p)$ induces a homeomorphism $V_0(R/(p))\rightarrow 
V_0(R)$.  Thus by Proposition~10.1, for 
i) it suffices to show that the map from 
$H^\spot(X)/(p)$ to $H^\spot(X;\fp)$ is an F-isomorphism.  The kernel
of this map is trivial, and it was shown during the proof of
Proposition~10.2 that if $p^j$ annihilates the torsion in $H^*(X)$,
then for any $x\in H^\spot(X;\fp)$, $x^{p^{j+1}}$ is in the image of
$H^\spot(X)/(p)$.  Part ii) also follows easily from Proposition~10.1,
because the kernel of the inclusion of $H^\spot(X;\ftwo)$ in
$H^*(X;\ftwo)$ is of course trivial, and the square of any element is
in the image.  

For iii), we first consider the case $i=0$.  Let $f$ stand for the map
from $H^*(X)\otimes H^*(X')$ to $H^*(X\times X')$.  Since $X$ and $X'$
are finite type CW-complexes, there is a K\"unneth exact sequence 
$$0\mapright{} H^*(X)\otimes H^*(X')\mapright{f}
H^*(X\times X')\mapright{} \Tor^{*-1}(H^*(X),H^*(X'))
\mapright{} 0\eqno{(6)}$$
which splits, and the $\Tor$-term consists of torsion elements of
bounded exponent.

Consider the commutative diagram given below.   
Each of the maps is injective (for $f/(p)$ this follows from the fact
that the K\"unneth sequence~(6) is split).  The maps $\pi$ and $\pi\otimes\pi$
are F-isomorphisms by the proof of i), and the map labelled \lq$\rm inc.$'
is an F-isomorphism because it is injective and the $p$th power of any
element is in its image.  Now $f/(p)$ is injective, and every other
map in the diagram is an F-isomorphism.  It follows that $f/(p)$ is
also an F-isomorphism.

\midinsert
$$\matrix{ 
(H^\spot(X)\otimes H^\spot(X'))/(p)&\mapright{\cong} &
H^\spot(X)/(p)\otimes H^\spot(X')/(p)\cr
\mapdown{f/(p)}&&\mapdown{\pi\otimes\pi}\cr
H^\spot(X\times X')/(p)&&
H^\spot(X;\fp)\otimes H^\spot(X';\fp)\cr
\mapdown{\pi}&&\mapdown{\rm inc.}\cr
H^\spot(X\times X';\fp)&\mapright{\cong}&
(H^*(X;\fp)\otimes H^*(X';\fp))^{\rm even}\cr}$$
\endinsert

Since $V_i(R)$ is a subvariety of $V_0(R)$, it follows that for
general $i$ the map 
$$f:V_i(H^\spot(X\times X'))\rightarrow V_i(H^\spot(X)\otimes
H^\spot(X'))$$
is a homeomorphism onto its image.  To show that this map is
surjective, it suffices to show that if $y\in H^\spot(X\times X')$
satisfies $p^iy=0$, then there exists $x\in
H^\spot(X)\otimes H^\spot(X')$ with $p^ix=0$ and such that $f(x)=y^N$
for some~$N$.  (This is because $\phi\in V_0(R)\setminus V_i(R)$ if and only
if there exists $r$ with $p^ir=0$ and $\phi(r)\neq 0$.)  From 
the fact that $f/(p)$ is an F-isomorphism, we obtain $x'\in
H^\spot(X)\otimes H^\spot(X')$ and $y'\in H^\spot(X\times X')$ such
that $f(x')=y^{p^n}+py'$.  Now pick $r$ sufficiently large that $p^r$
annihilates the torsion in $H^*(X)$ and $H^*(X')$.  Then $p^r$
annihilates the $\Tor$-term in~(6), so $f(x^{\prime p^r})-y^{p^{n+r}}$
is in the image of $f$.  Thus there exists $x''\in H^*(X)\otimes 
H^*(X')$ such that $f(x'') = y^{p^{n+r}}$.  Express $x''= x_\eve+x_\odd$,
where $x_\eve\in H^\spot\otimes H^\spot$ and $x_\odd\in H^{\rm odd}\otimes 
H^{\rm odd}$.  Now $p^ix_\eve= p^ix_\odd=0$, and either $p$ is odd in which
case $x_\odd^2=0$, and 
$$x^{\prime\prime p^i}= x_\eve^{p^i} + p^ix_\eve^{p^i-1} + x_\odd^2(\cdots)  
= x_\eve^{p^i},$$
or $p=2$, in which case $x_\odd^2\in H^\spot(X)\otimes H^\spot(X')$,
$2x_\odd^2=0$, and 
$$x^{\prime\prime 2^i} = x_\eve^{2^i} + 2^i x_\eve^{2^i-1}x_\odd +
2x_\odd^2(\cdots) + x_\odd^{2^i} = x_\eve^{2^i} + x_\odd^{2^i}.$$
In either case, $x= x^{p^i}$ is an element of $H^\spot(X)\otimes
H^\spot(X')$ such that $f(x) = y^{p^N}$ for some~$N$, as required.  
\qed

\proclaim Proposition 10.4. Let $R$ be a finitely generated commutative
$\plocal$-algebra, with an action of a finite group $G$, and write
$R^G$ for the $G$-fixed points in $R$.  Then for
each $i$, the natural map gives rise to a homeomorphism 
$V_i(R)/G\rightarrow V_i(R^G)$.  

\proof For $i=0$, this is a standard result, see for example chapter~5
of [\atimac].  It is also a special case of Lemma~8.11 of [\quib]
(view the group $G$ as a category with one object and the $G$-action
on  $R$ as a functor from $G$ to $\plocal$-algebras).  Given the case
$i=0$, the general case follows provided that $V_i(R)$ maps
surjectively to $V_i(R^G)$.   Let $\phi$ be any element of $V(R)$
whose restriction to $R^G$ lies in $V_i(R^G)$, and let $r$ be an
element of $R$ such that $p^ir=0$.  Then $Q_r(r)=0$, where $Q_r(X)$ is
the monic polynomial 
$$Q_r(X) = \prod_{g\in G} (X-g.r).$$
The coefficients in this polynomial lie in $R^G$, and (apart from 1,
the leading coefficient) are annihilated by $p^i$.  Thus 
$$\phi(r^{|G|}) = \phi(f_r(r)) = \phi(0) = 0,$$
and so $\phi\in V_i(R)$.  
\qed

\proclaim Proposition 10.5.  Let $R$ be a finitely generated
commutative $\plocal$-algebra, let $C_p$ act on $S=R^{\otimes p}$ by
permuting the factors, and let $\mu : S\rightarrow R$ be the
multiplication homomorphism.  Define $U\subseteq V(S)$ by 
$$U=\{ \phi\in V(S) \mid \forall x\in S, 
\phi\Bigl(\sum_{g\in C_p} g.x\Bigr) =0\}.$$ 
Then $U= \mu^*(V(R))$, and $U$ maps injectively to $V(S)/C_p=
V(S^{C_p})$.  

\proof First note that $\mu:S\rightarrow R$ is $C_p$-equivariant for
the given action on $S$ and the trivial action on $R$.  Now if 
$\phi\in \mu^*(V(R))$, or equivalently, $\phi= \psi\circ\mu$ for some 
$\psi : R\rightarrow k$, then for any $x\in S$, 
$$\phi\Bigl(\sum_{g\in C_p} g.x\Bigr)= 
\psi\Bigl(\sum_{g\in C_p} g.\mu(x)\Bigr) = \psi(p.\mu(x))= 0.$$
Thus $\mu^*(V(R))\subseteq U$.  

For the converse, given $r\in R$, let $x_1= r\otimes 1\otimes
\cdots\otimes 1$, and let $x_2,\ldots,x_p$ be the images of $x_1$
under the $C_p$-action.  It suffices to show that if $\phi\in U$, then
for any $r$, $\phi(x_1)=\phi(x_2)=\cdots =\phi(x_p)$.  
Let $\sigma_i=\sigma_i(x_1,\ldots,x_p)$ be the $i$th elementary
symmetric function in $x_1,\ldots,x_p$.  For $1\leq i\leq p-1$, the
stabilizer in $\Sigma_p$ of a subset of size $i$ has order coprime to
$p$, and hence there exists $y_i$ such that 
$$\sigma_i= \sum_{g\in C_p}g.y_i.$$
Now if $\phi\in U$, then $\phi(\sigma_i)=0$ for $1\leq i\leq p-1$, and
hence each $\phi(x_j)$ is a root of the equation
$X^p-\phi(\sigma_p)=0$, which has all of its roots equal.  

The last claim follows since $\mu^*(V(R))$ consists of points of
$V(S)$ fixed by the $C_p$-action.  \qed

\proclaim Theorem 10.6.  If $X$ is a finite type CW-complex such that
$H^*(X)$ is finitely generated as a $\plocal$-algebra, then so is
$X^p\times_{C_p}EC_p$, and $V_i(H^\spot(X^p\times_{C_p}EC_p))$ may be
described as follows:  
\pra 
Let $\gamma^*$, $\delta^*$, be the ring homomorphisms 
$$H^*(X\times BC_p)\mapright{\gamma^*}H^*(X)\qquad 
(H^*(X)^{\otimes p})^{C_p}\mapright{\delta^*}H^*(X)$$
induced by the projection $X\times EC_p\rightarrow X\times BC_p$ and
the diagonal map $X\times EC_p\rightarrow X^p\times EC_p$.  Then 
$$V_0(H^\spot(X^p\times_{C_p}EC_p))\cong
\lim_{\rightarrow}(V_0(H^\spot(X\times BC_p))\mapleft{\gamma} 
V_0(H^\spot(X))\mapright{\delta}V_0(H^\spot(X)^{\otimes p})/{C_p}),$$
and for $i>0$, 
$$V_i(H^\spot(X^p\times_{C_p}EC_p))\cong 
\lim_{\rightarrow}(V_{i-1}(H^\spot(X))\mapleft{\rm Id} 
V_i(H^\spot(X))\mapright{\delta}V_i(H^\spot(X)^{\otimes p})/{C_p}).$$

\remark Roughly one could write 
$$\eqalign{ 
V_0(H^\spot(X^p\times_{C_p}EC_p)) &= V_0(H^\spot(X\times BC_p)) \cup 
V_0(H^\spot(X)^{\otimes p})/{C_p},\cr
V_i(H^\spot(X^p\times_{C_p}EC_p)) &= V_{i-1}(H^\spot(X)) \cup 
V_i(H^\spot(X)^{\otimes p})/{C_p}.\cr}$$

\proof Let $\alpha^*$ and $\beta^*$ be the homomorphisms defined in
Section~6, so that there is a commutative diagram:  
$$\matrix{
H^*(X^p\times_{C_p}EC_p)&\mapright{\alpha^*}&H^*(X^p)^{C_p}\cr
\mapdown{\beta^*}&&\mapdown{\delta^*}\cr 
H^*(X\times BC_p)&\mapright{\gamma^*}&H^*(X).\cr
}$$
Use the same notation for the corresponding homomorphisms of mod-$p$
cohomology.  

Obviously, $X$ is of finite type if and only if
$X^p\times_{C_p}EC_p$ is of finite type.  In this case, Nakaoka's
theorem given above as Theorem~2.1 (together with the remarks that
follow it), implies that $H^*(X^p\times_{C_p}EC_p;\fp)$ is isomorphic to
the graded tensor product 
$$H^*(X^p;\fp)^{C_p}\otimes H^*(BC_p;\fp)$$
modulo the ideal generated by elements of the form 
$$\sum_{g\in C_p} g^*(x)\otimes y,\qquad x\in H^*(X^p;\fp),\quad
y \in H^i(BC_p;\fp), \quad i>0.$$
Let $\eta^*$ stand for the projection map 
$$\eta^*:H^*(X^p;\fp)^{C_p}\otimes H^*(BC_p;\fp)\mapright{}
H^*(X^p\times_{C_p}EC_p;\fp).$$
Note that there is no analogue of $\eta^*$ for $\plocal$-cohomology.  
Note also that there is a commutative diagram:  
$$\xymatrix{H^*(X^p;\fp)^{C_p}\otimes H^*(BC_p;\fp)\ar[r]^{\eta^*}
\ar[rd]_{\delta^*\otimes 1} &
H^*(X^p\times_{C_p}EC_p;\fp)\ar[d]^{\beta^*}\\
&H^*(X;\fp)\otimes H^*(BC_p;\fp).}$$

It follows that $H^*(X^p\times_{C_p}EC_p;\fp)$ is finitely generated if
$H^*(X;\fp)$ is.  If $p^i$ annihilates torsion in $H^*(X)$, then
either by the results of Section~4 or a transfer argument, $p^{i+1}$
annihilates the torsion in $H^*(X^p\times_{C_p}EC_p)$.  Hence by
Proposition~10.2, if $H^*(X)$ is finitely generated, then so is 
$H^*(X^p\times_{C_p}EC_p)$.  

For the case $i=0$, it is convenient to work with mod-$p$ cohomology,
and deduce the result for $\plocal$-cohomology from Proposition~10.3.  
Let $\phi$ be an element of $V_0(H^*(X^p\times_{C_p}EC_p;\fp))$, and
let $y$ be the image in $H^2(X^p\times_{C_p}EC_p;\fp)$ of a generator
for $H^2(BC_p;\fp)$.  If $y'$ is any element of
$H^*(X^p\times_{C_p}EC_p;\fp)$ such that $\alpha^*(y')=0$, 
then $y^{\prime 2}$ is in the ideal generated by $y$.  (To see this,
recall that $\alpha^*$ is the edge homomorphism in the Cartan-Leray
spectral sequence for $X^p\times_{C_p}EC_p$, and that the ideal
generated by $y$ maps on to $E_2^{i,j}$ for $i\geq 2$.)  Thus if
$\phi(y)= 0$, then $\phi$ factors through $H^\spot(X^p;\fp)^{C_p}$, 
or in other words $\phi \in \alpha(V(H^\spot(X^p;\fp)^{C_p}))$, 
and this factorisation is of course unique.  

On the other hand, if $\phi(y)\neq 0$, then for all $x\in
H^\spot(X^p;\fp)$, 
$$\phi(\sum_{g\in C_p}g^*(x))= 0.$$
In particular, Proposition~10.5 in the case $R=H^\spot(X;\fp)$ shows
that there exists $\phi'\in V(H^\spot(X;\fp)\otimes
H^\spot(BC_p;\fp))$ such that $\phi\circ\eta^*=
\phi'\circ(\delta^*\otimes 1)$.  Since $\eta^*$ is surjective, it
follows that $\phi = \phi'\circ \beta^*$.  

Thus if $\phi(y)\neq 0$ then $\phi$ factors
through the image of $\beta^*$.  If $R$ is the image of $H^\spot(X;\fp)$
under the map $x\mapsto x^p$, then the image of $\beta^*$ 
contains $R\otimes H^\spot(BC_p;\fp)$, as this is the subring
generated by the image of $H^\spot(BC_p;\fp)$ and elements of the form
$x \wr 1$.  However, any homomorphism from $R\otimes
H^\spot(BC_p;\fp)$ to $k$ extends uniquely to $H^\spot(X\times
BC_p;\fp)$ by Proposition~10.1, and so $\phi$ factors uniquely through
$H^\spot(X\times BC_p;\fp)$.  

Finally, if $\phi$ factors through both $H^\spot(X^p;\fp)^{C_p}$ and 
$H^\spot(X\times BC_p;\fp)$, then $\phi(y)=0$, and so $\phi$ factors
uniquely through $H^\spot(X;\fp)$.  This completes the claim in the
case $i=0$.  

For the case when $i>0$, we consider the information given by the
Cartan-Leray spectral sequence with $\plocal$-coefficients.  The 
results of Section~5 imply that any element $x$ of
$H^*(X^p\times_{C_p}EC_p)$ yielding an element of $E_\infty^{m,n}$ for
$m>0$, has order $p$.  Hence if $\phi$ is an element of
$V_1(H^\spot(X^p\times_{C_p}EC_p))$, $\phi$ factors through
$H^\spot(X^p)^{C_p}$.  Thus for each $i>0$, 
$$\alpha(V(H^\spot(X^p)^{C_p}))\supseteq 
V_i(H^\spot(X^p\times_{C_p}EC_p))\supseteq 
\alpha(V_i(H^\spot(X^p)^{C_p})).$$
For simplicity we shall restrict attention to
$H^{2p*}$, the subring of $H^\spot$ consisting of elements in degrees
divisible by $2p$.  (The inclusion of $H^{2p*}$ in $H^\spot$ induces
a homeomorphism from $V_i(H^\spot)$ to $V_i(H^{2p*})$ for each $i$.)  

As in previous sections we fix a decomposition of $H^*(X)$ as a direct
sum of cyclic groups, and use this to split the $E_\infty$-page of the
Cartan-Leray spectral sequence as a direct sum of pieces as described
in Section~5.  Refine the resulting decomposition of
$H^*(X^p\times_{C_p}EC_p)$ to a splitting into cyclic summands.  
If $x$ generates a cyclic summand of
$H^{2p*}(X^p\times_{C_p}EC_p)$ of order $p^i$, there are only two 
possibilities:  Either $x$ is contained in a summand of the spectral
sequence concentrated in $E_*^{0,*}$, and its image in
$H^{2p*}(X^p)^{C_p}$ generates a cyclic summand of order $p^i$, or 
$x$ is contained in a summand of the spectral sequence of the type 
discussed in Theorems 5.1~and~5.1$'$, and its image in
$H^{2p*}(X^p)^{C_p}$ generates a cyclic summand of order $p^{i-1}$.
In the first case, there exists $x'\in H^{2p*}(X^p)$ such that 
$$x=\sum_{g\in C_p} g.x',$$
and so (as in the proof of Proposition~10.5) $\delta^*(x)=0$.  
It follows that 
$$V_i(H^{2p*}(X^p\times_{C_p}EC_p))\supseteq 
\alpha(V_i(H^{2p*}(X^p)^{C_p}))\cup 
\alpha\circ\delta(V_{i-1}(H^{2p*}(X)), $$
and it remains only to prove the opposite inequality.  

If $\phi\in V(H^{2p*}(X))\setminus V_{i-1}$, then there exists $x\in
H^{2p*}(X)$ such that $p^{i-1}x=0$ and $\phi(x)\neq 0$.  In this case,
$y = x\wr 1$ is an element of $H^{2p*}(X^p\times_{C_p}EC_p)$ such that 
$p^i y = 0$ and $\phi\circ\delta^*\circ\alpha^*(y) \neq 0$.  Thus
$$\alpha\circ\delta(\phi)=\phi\circ\delta^*\circ\alpha^*
\notin V_i(H^{2p*}(X^p\times_{C_p}EC_p)).$$  
If $$\phi \in V(H^{2p*}(X^p))\setminus (\delta(V(H^{2p*}(X)))\cup 
V_i(H^{2p*}(X^p))),$$ then there exists $x,x'\in H^{2p*}(X)^{\otimes p}$
such that $p^i.x=0$ but $\phi(x)\neq 0$, and (by Proposition~10.5) 
$\phi\bigl(\sum_{g\in C_p}g.x'\bigr)\neq 0$.  It suffices to find
$y\in H^{2p*}(X)^{\otimes p}$ such that $p^i.y =0$ and 
$\phi\bigl(\sum_{g\in C_p}g.y\bigr)\neq 0$.  By choosing $y_j$ to be a
combination of various products of elements of the form $g.x$, we may
ensure that for $1\leq j\leq p-1$, 
$$\sum_{g\in C_p}g.y_j = \sigma_j(g_1.x,\ldots,g_p.x),$$
where $g_1.x,\ldots,g_p.x$ are the images of $x$
under the $C_p$-action, and $\sigma_j$ stands for the $j$th elementary
symmetric function.  Each $y_j$ satisfies $p^i.y_j=0$, so if there
exists $j$ such that 
$$\phi\Bigl(\sum_{g\in C_p}g.y_j\Bigr) = 
\phi\Bigl(\sigma_j(g_1.x,\ldots,g_p.x)\Bigr)\neq 0, $$
then $y=y_j$ will do.  Otherwise, it follows (as in the proof of
Proposition~10.5) that for each $g\in C_p$, $\phi(g.x)= \phi(x)$.  
In this case, $y=xx'$ has the required properties.  
\qed

\beginsection 11. CARLSON'S $W_i(G)$  

As in the last section we take $\plocal$-coefficients for cohomology
unless otherwise stated, and fix an algebraically closed field $k$ of
characteristic $p$, upon which the definition of $V(-)$ depends.  
For $G$ a finite group, note that $EG$ may
be taken to be of finite type.  For $G$ a finite group and 
$Y$ a finite $G$-CW-complex, Quillen showed
that the equivariant cohomology ring $H^*(Y\times_GEG)$ is finitely
generated [\quib].  For such $G$ and $Y$, define varieties $W_i(G,Y)$,
$W_i(G)$ by 
$$W_i(G,Y)= V_i(H^\spot(Y\times_GEG)),\qquad W_i(G) = W_i(G,\{*\}),$$
where $V_i(-)$ is as defined in Section~10.  
Note that $W_i(G,Y)$ is a covariant functor of the pair $(G,Y)$, 
and that if $H$ is a subgroup of $G$, then the finiteness of $H^*(H)$
as an $H^*(G)$-module implies that the fibres of the map
$W_i(H)\rightarrow W_i(G)$ are finite.  
In [\quib], Quillen gave a complete description of $W(G,Y) =
W_0(G,Y)$.  The subvarieties $W_i(G)$ for $i>0$ were introduced by 
Carlson in [\carl].  Very little is known about $W_i(G)$ and
$W_i(G,Y)$ in general.  The results of Section~10 of course have
Corollaries concerning $W_i(G,Y)$.  For example, if $Y'$ is a finite
$H$-CW-complex, then by Proposition~10.3~iii), 
$$W_i(G\times H,Y\times Y') \cong W_i(G,Y)\times W_i(H,Y').$$
For most of this section we use examples constructed as wreath
products to shed some light on questions raised by Carlson concerning
$W_i(G)$.  The lower bound we give for the size of $W_1(G)$ in
Theorem~11.9 is independent of the rest of the paper however.  
The following three statements are corollaries of Theorem~10.6.

\proclaim Corollary 11.1.  Let $G$ be a finite group and $Y$ a finite
$G$-CW-complex, and let $G\wreath C_p$ act on $Y^p$, where $G^p$ acts
component-wise and $C_p$ by permuting the factors.  Let $\gamma$ be
the inclusion of the $G$-space $Y$ as the diagonal in $Y^p$, which is
equivariant for the diagonal map from $G$ to $G^p$, let $C_p$ act
trivially on $Y$, and let $\delta$ be the identity map on $Y$ viewed
as an equivariant map from the $G$-space $Y$ to the $G\times
C_p$-space $Y$.  Then 
$$W_0(G\wreath C_p,Y^p) = 
\lim_{\rightarrow}(W_0(G^p,Y^p)/C_p\mapleft{\gamma} W_0(G,Y)
\mapright{\delta}W_0(G\times C_p,Y)),$$
and for $i>0$, 
$$W_i(G\wreath C_p,Y^p) = 
\lim_{\rightarrow}(W_i(G^p,Y^p)/C_p\mapleft{\gamma} W_i(G,Y)
\mapright{\rm Id}W_{i-1}(G,Y)).$$

\proof This is just Theorem 10.6 in the case $X= Y\times_GEG$.  
\qed

\proclaim Corollary 11.2.  With notation as in Corollary~11.1, 
and $i>0$, 
$$\eqalign{\dim W_0(G\wreath C_p,Y^p) 
&= \max\{ p.\dim W_0(G,Y),\, 1\},\cr
\dim W_i(G\wreath C_p,Y^p)
&= \max\{p.\dim W_i(G,Y),\, \dim W_{i-1}(G,Y)\}.\cr}$$
\qed

\proclaim Proposition 11.3.  Let $m= \sum_{j}m_jp^j$, where $0\leq
m_j<p$.  Then 
$$\dim W_i(\Sigma_m)= \sum_{j\geq i+1}m_jp^{j-i-1}.$$
In particular, 
$$\dim W_i(\Sigma_{p^n}) = \cases{0&if $i \geq n$,\cr
p^{n-i-1}&otherwise.\cr}$$

\proof If $G$ is a finite group with Sylow $p$-subgroup $P$, then the
map $W_i(P)\rightarrow W_i(G)$ has finite fibres (by the Evens-Venkov
theorem) and is surjective since the kernel of the transfer is an
ideal forming a complement to the image of $H^*(G)$ in $H^*(P)$.  
The Sylow $p$-subgroup of $\Sigma_m$ is isomorphic to 
$$ (P_1)^{m_1}\times(P_2)^{m_2}\times \cdots,$$
where $P_n$ is the Sylow $p$-subgroup of $\Sigma_{p^n}$.  Thus it
suffices to prove the assertion for $P_n$.  This follows from
Corollary~11.2 by induction, since $P_n\cong P_{n-1}\wreath C_p$.  
\qed

Carlson also introduced the cohomological exponential
invariant, $che(G)$ of a finite group $G$, defined by 
$$\eqalign{
che(G) = \min\{O(x_1)&\cdots O(x_r) \mid H^*(G;\Bbb Z) 
\hbox{ is a finitely generated}\cr
&\hbox{module for the subring generated by the $x_i$}\}.}$$
The theorem that we quote as Theorem~3.4 is
equivalent to the statement that $|G|$ divides $che(G)$.  Let
$\nu_p(che(G))$ stand for the number of times that $p$ divides
$che(G)$.  Carlson showed that 
$$\nu_p(che(G))= \sum_{i\geq 0} \dim W_i(G)\eqno{(7)},$$
and gave an example of a group $G$ of order $p^5$ with $che(G)=p^6$.
Using Corollary~11.2 and the description of (the $p$-part of) $che(G)$
given by~(7), it is easy to find groups for which $|G|<che(G)$.  One
such example is contained in Proposition~11.4.  We give another
different example in Proposition~11.10.  

\proclaim Proposition 11.4. Let $G= (C_p^n)\wreath C_p$, so that $G$
is a split extension with kernel $(C_p)^{pn}$, quotient $C_p$, and the
quotient acts freely on an $\fp$-basis for the kernel.  Then 
$$che(G)=p^{pn+n},\quad\hbox{whereas}\quad |G|=p^{pn+1}.$$
\qed

The generalized Frattini subgroups $\Phi_n(G)$, 
defined in Section~8 (just before
Proposition~8.4) may be used to give an upper bound for $W_n(G)$.
Before stating our result in Theorem~11.6, we recall a theorem of 
Carlson from [\carl, section~4].

\proclaim Theorem 11.5. Let $G$ be a finite group and $H$ a normal
subgroup of $G$.  The inverse image of\/ $0$ under the map
$W(G)\rightarrow W(G/H)$ is equal to the image of $W(H)$.  
\qed

\remark The element $0\in W(G)$ is the homomorphism $H^*(G)\rightarrow
k$ which sends all elements in positive degree to zero.  

\proclaim Theorem 11.6. Let $G$ be a $p$-group.  Then $W_n(G)$ is
contained in the image of\/ $W(\Phi_n(G))$.  

\proof The group $H^i(G/\Phi_n(G))$ has exponent dividing $p^n$ for
all but finitely many $i$ by Proposition~8.4.  Thus if $\phi\in
W_n(G)$, then $\phi$ maps to $0$ in $W(G/\Phi_n(G))$.  The claim
follows from Theorem~11.5.  
\qed

Carlson asked if, for $G$ a $p$-group, $W_n(G)$ always contains the
image of $W(\Phi^n(G))$, where $\Phi^n(G)=\Phi(\Phi^{n-1}(G))$ is the
$n$th iterated Frattini subgroup of $G$.  
Proposition~11.7 shows that this
is consistent with the upper bound of Theorem~11.6.  However, in
Proposition~11.8 we construct wreath products $G$ such that $W_n(G)$
does not contain the image of $W(\Phi^n(G))$, answering Carlson's
question negatively.  

\proclaim Proposition 11.7. Let $G$ be a $p$-group.  Then
$\Phi_n(G)\geq \Phi^n(G)$.  

\proof Given $H$, a subgroup of $G$ of index $p^r\leq p^n$, let 
$$H=H_r<H_{r-1}<\cdots<H_1<H_0=G$$
be a chain of subgroups of $G$ such that $|H_i:H_{i+1}|=p$.  Then for
each $i$, $H_{i+1}\geq \Phi(H_i)$, and so by induction 
$H_r\geq \Phi^r(G)\geq \Phi^n(G)$.  
\qed

\proclaim Proposition 11.8. Let $S$ be an elementary abelian $p$-group of
rank $r$, acting freely transitively on $\Omega$, and let $G$ be the
wreath product $C_{p^n}\wreath S$ for some $n>0$.  Then 
$$\dim W(\Phi^n(G))= p^r-1,\quad\hbox{whereas}\quad 
\dim W_n(G)\leq p^{r-1}.$$
In particular, $W_n(G)$ cannot contain the image of $W(\Phi^n(G))$
unless $p=2$ and $r=1$.  

\proof Recall that for $H$ a subgroup of $G$, the map $W(H)$ to $W(G)$
is finite, and so preserves dimensions.  Recall also that for $G$ a
$p$-group, $\Phi(G)$ can be defined to be the minimal normal subgroup
$H$ of $G$ such that $G/H$ is elementary abelian.  Now $G$ as in the
statement can be generated by $r+1$ elements, $r$ of which map to a
generating set for $S$, and one element of the form 
$$(g,1,\ldots,1)\in (C_{p^n})^\Omega \leq C_{p^n}\wreath S=G.$$
A homomorphism from $G$ onto $(C_p)^{r+1}$ may be constructed, which
shows that $G$ cannot be generated by fewer than $r+1$ elements.  It
follows that $\Phi(G)$ has index $p^{r+1}$, and is an index $p$ 
subgroup of $(C_{p^n})^\Omega$.  Thus 
$$\Phi(G)\cong (C_{p^n})^{p^r-1}\times C_{p^{n-1}},
\qquad
\Phi^n(G)\cong (C_p)^{p^r-1}.$$
This proves the first claim.  For the second claim, note that $G$ is a
subgroup of the wreath product 
$C_{p^n}\wreath \Sigma_{p^r}\leq \Sigma_{p^{n+r}}$, 
and so $\dim W_n(G)$ is bounded by $\dim W_n(\Sigma_{p^{n+r}})$, which
equals $p^{r-1}$ by Proposition~11.3.  
\qed

The only general lower bound that I have been able to find is the
following bound for $W_1(G)$, which is related to Proposition~1 
of [\ian].  

\proclaim Theorem 11.9. Let $G$ be a $p$-group, and let $H=Z(G)\cap
\Phi(G)$, the intersection of the centre and the Frattini subgroup 
of $G$.  Then $W_1(G)$ contains the image of $W(H)$.  

\proof $H$ is central, so is abelian.  Let $K$ be the elements of $H$
of order dividing $p$.  Then $K$ is a subgroup of $H$ and $W(K)$ maps
homeomorphically onto $W(H)$.  It will therefore suffice to show that 
any element of $H^*(G)$ of order $p$ has trivial image in $H^*(K)$.  

For this, we claim that the image of $H^*(G;\fp)$ in $H^*(K;\fp)$ is
contained in the image of $H^*(K)$.  Assuming this for now, let $x\in
H^*(G)$ have order $p$, and recall the Bockstein long exact sequence 
$$\cdots H^*(G)\mapright{\times p}H^*(G)\mapright{}
H^*(G;\fp)\mapright{\delta}H^{*+1}(G)\mapright{\times p}
H^{*+1}(G)\cdots.$$
Since $px=0$, there exists $y$ such that $\delta(y)=x$.  Now 
$\delta(\res^G_K(y))= \res^G_K(x)$.  By the claim, the image of
$H^*(G;\fp)$ in $H^*(K;\fp)$ is contained in the kernel of $\delta$,
and so $\res^G_K(x)=0$ as required.  

It remains to prove the claim made in the last paragraph.  For this we
consider the spectral sequence with $\fp$-coefficients for the central
extension 
$$K\mapright{}G\mapright{}G/K.$$
The $E_2$-page is isomorphic to $H^*(K;\fp)\otimes H^*(G/K;\fp)$.  For
clarity, suppose that $p$ is odd (the case $p=2$ is similar but not
identical).  In this case, 
$$H^*(K;\fp)\cong \Lambda[K^\#]\otimes \fp[K],$$
the tensor product of the exterior algebra on $K^\#=\Hom(K,\fp)$,
generated in degree 1, with the algebra of polynomial functions on $K$
(where the monomials have degree 2).  The Bockstein
$H^*(K;\fp)\rightarrow H^{*+1}(K;\fp)$ maps $H^1$ isomorphically to
the degree 2 piece of $\fp[K]$.  

Since $K$ is contained in $\Phi(G)$, every element of $H^1(G;\fp)\cong
\Hom(G,\fp)$ restricts trivially to $K$.  It follows that 
$$d_2:E_2^{0,1}\cong K^\#\mapright{}E_2^{2,0}\cong H^2(G/K;\fp)$$
must be injective.  Note also that the generators for $\fp[K]$ are
cycles for $d_2$ by the Serre transgression theorem.  It follows that 
$E_3^{0,*}$, the cycles for $d_2$ in $E_2^{0,*}$, is isomorphic to
$\fp[K]$.  This subring of $H^*(K;\fp)$ is contained in the image of
$H^*(K)$.  The image of $H^*(G;\fp)$ in $H^*(K;\fp)$ is isomorphic to 
$E_\infty^{0,*}$, so is a subring of the image of $H^*(K)$ as
required.  
\qed

\remark Note that Theorems 11.6 and 11.9 together describe $W_1(G)$
completely for any $p$-group $G$ such that $\Phi(G)$ is central.  
This includes the two examples discussed by Carlson in [\carl].  As
a further example we make the following statement, whose proof we leave as
an exercise.  

\proclaim Proposition 11.10.  Let $G$ be the group with presentation 
$$G=\langle a_1,\ldots,a_n \mid a_i^p=1, [a_i,a_j] \hbox{ is
central}\rangle.$$
Then $G$ is a finite $p$-group, with 
$$\dim W_0(G) = {n \choose 2} +1, \qquad \dim W_1(G) = {n\choose 2},$$
$$\nu_p(|G|)={n+1\choose 2}, \qquad \nu_p(che(G))= n^2-n+1.$$
\qed

\remark The author knows of no example in which $W_i(G)$ is not equal
to the image of $W(H)$ for some subgroup $H$ of $G$.  One might
conjecture that this is always the case, or more weakly one might
conjecture that for any $G$, $W_i(G)$ is equal to a union of images of
$W(E)$ for some collection of elementary abelian subgroups $E$ of
$G$.  By section~12 of [\quib], this weaker conjecture is equivalent
to the following:  Let $J_i$ be the ideal of $H^*(G;\fp)$ generated by
the image of $I_i$, the ideal in $H^*(G)$ used to define $W_i(G)$.  
Then the radical of $J_i$ is closed under the action of Steenrod's
reduced powers ${\cal P}^i$.  

\beginsection REFERENCES

\frenchspacing
\def\book#1/#2/#3/#4/#5/{\item{#1} #2, {\it #3,} #4, {\oldstyle #5}.
\par\smallskip}
\def\sebook#1/#2/#3/#4/#5/#6/#7/{\item{#1} #2, {\it #3,} #6, {\oldstyle #7}.
\par\smallskip}
\def\paper#1/#2/#3/#4/#5/(#6) #7--#8/{\item{#1} #2, #3, {\it #4,} {\bf #5}
({\oldstyle#6}) {\oldstyle #7}--{\oldstyle#8}.\par\smallskip}
\def\prepaper#1/#2/#3/#4/#5/#6/{\item{#1} #2, #3, {\it #4} {\bf #5} {#6}.
\par\smallskip}

\def\book#1/#2/#3/#4/#5/{\item{#1.} #2, \lq\lq{#3},'' #4, #5.
\par\smallskip}
\def\sebook#1/#2/#3/#4/#5/#6/#7/{\item{#1.} #2, \lq\lq{#3},'' 
#6, #7. \par\smallskip}
\def\paper#1/#2/#3/#4/#5/(#6) #7--#8/{\item{#1.} #2, #3, {\it #4} {\bf #5}
(#6), #7--#8. \par\smallskip}
  
\paper \ade/A. Adem/Cohomological exponents of $\Bbb ZG$-lattices/J. Pure
Appl. Alg./58/(1989) 1--5/

\paper \adca/A. Adem and J. F. Carlson/Discrete groups with large
exponents in cohomology/J. Pure Appl. Alg./66/(1990) 111--120/

\sebook \ademil/A. Adem and R. J. Milgram/Cohomology of Finite
Groups/Grundlehren der mathematischen Wissenschaften/309/
Springer Verlag/1994/

\book \atimac/M. F. Atiyah and I. G. Macdonald/Introduction to
Commutative Algebra/Addison-Wesley/1969/ 

\sebook \ben/D. J. Benson/Representations and Cohomology Vols. I and
II/Cambridge Studies in Advanced Mathematics/30 {\rm and} 
31/Cambridge Univ. Press/1991/

\sebook \beninv/D. J. Benson/Polynomial Invariants of Finite Groups/
London Math. Soc. Lecture Notes/190/Cambridge Univ. Press/1993/

\item{\bencarl.} D. J. Benson and J. F. Carlson, The cohomology of
extraspecial groups, {\it Bull. London Math. Soc.} {\bf 24}
({1992}), {209}--{235} and erratum, 
{\it Bull. London Math. Soc.} {\bf 25} ({1993}), 498.  
\par\smallskip

%\paper \bencarl/D. J. Benson and J. F. Carlson/The cohomology of
%extraspecial groups/Bull. London Math. Soc./24/(1992) 209--235/

\sebook \bro/K. S. Brown/Cohomology of Groups/Graduate Texts in
Mathematics/87/Springer Verlag/1982/

\book \brobui/K. S. Brown/Buildings/Springer Verlag/1989/

\sebook \lmm/R. R. Brunner, J. P. May, J. E. McClure and M.
Steinberger/${\rm H}_\infty$ Ring Spectra and their Applications/Lecture
Notes In Mathematics/1176/Springer Verlag/1986/

\paper \carl/J. F. Carlson/Exponents of modules and maps/Invent.
Math./95/(1989) 13--24/

\book \caei/H. Cartan and S. Eilenberg/Homological Algebra/Princeton
Univ. Press/1956/

\paper \eva/L. Evens/On the Chern classes of representations of finite
groups/Trans. Amer. Math. Soc./115/(1965) 180--193/

\sebook \evb/L. Evens/The Cohomology of Groups/Oxford Mathematical
Monographs//Clarendon Press/1991/

\paper \evk/L. Evens and D. S. Kahn/Chern classes of certain
representations of symmetric groups/Trans. Amer. Math. Soc./245/(1978)
309--330/ 

\paper \hako/M. Harada and A. Kono/On the integral cohomology of
extraspecial $2$-groups/J. Pure Appl. Alg./44/(1987) 215--219/ 

\paper \hwh/H.-W. Henn/Classifying spaces with injective mod $p$
cohomology/Comment. Math. Helvetici/64/(1989) 200--206/ 

\paper \huea/J. Huebschmann/Cohomology of metacyclic groups/Trans.
Amer. Math. Soc./328/(1991) 1--72/

\item{\jrh.} J. R. Hunton, The complex oriented cohomology of extended
powers, preprint, 1995.  
\par\smallskip

\paper \kah/B. Kahn/Classes de Stiefel-Whitney de Formes quadratiques
et de repr\'esentations galoisiennes r\'eelles/Invent. Math./78/(1984)
223--256/ 

%\def\paper#1/#2/#3/#4/#5/(#6) #7--#8/{\item{#1} #2, #3, {\it #4,} {\bf #5}
%({\oldstyle#6}) {\oldstyle #7}--{\oldstyle#8}.\par\smallskip}
%\item{\ian} I. J. Leary, A bound on the exponent of the cohomology of
%$BC$-bundles, to appear in the proceedings of the Barcelona
%Conference on Algebraic Topology, 1994.  

%\paper \ian/I. J. Leary/A bound on the exponent of the cohomology of
%$BC$-bundles/in Algebraic topology: new trends in localization and
%periodicity edited by C. Broto, C. Casacuberta and G. Mislin

\item{\ian.} I. J. Leary, A bound on the exponent of the cohomology of
$BC$-bundles, {\it in} \lq\lq{Algebraic Topology: New Trends in Localization
and Periodicity}'' (eds. C. Broto, C. Casacuberta, G. Mislin), Progress in
Mathematics Vol. 136, Birkh\"auser, 1996.  
\par\smallskip

\book \oprob/J. van Mill and G. M. Reed (eds.)/Open Problems in
Topology/North Holland/1990/

\paper \nak/M. Nakaoka/Homology of the infinite symmetric group/Ann.
of Math./73/(1961) 229--257/

\paper \pal/F. P. Palermo/The cohomology ring of product complexes/Trans.
Amer. Math. Soc./86/(1957) 174--196/

\paper \qui/D. Quillen/The Adams conjecture/Topology/10/(1971) 67--80/

\item{\quib.} D. Quillen, The spectrum of an equivariant cohomology ring
I and II, {\it Ann. of Math.} {\bf 94} ({1971}), 549--572 and 573--602.  

\paper \quies/D. Quillen/The mod-$2$ cohomology rings of extra-special
$2$-groups and the spinor groups/Math. Ann./194/(1971) 197--212/  

\book \spa/E. Spanier/Algebraic Topology/McGraw-Hill/1966/

\sebook \ste/N. Steenrod and D. B. A. Epstein/Cohomology
operations/Ann. of Math. Studies/50/Princeton Univ. Press/1962/

\end